\documentclass[mnsc]{informs3} 

\OneAndAHalfSpacedXI

\usepackage{epstopdf}
\usepackage{natbib}
 \def\bibfont{\small}%
 \def\bibsep{\smallskipamount}%
 \def\bibhang{24pt}%
\TheoremsNumberedThrough     
\ECRepeatTheorems

\EquationsNumberedThrough    

\MANUSCRIPTNO{MS-0001-1922.65}

\begin{document}



\RUNTITLE{The user base dynamics of websites}

\TITLE{The user base dynamics of websites}

\ARTICLEAUTHORS{%
\AUTHOR{Kartik Ahuja}
\AFF{Electrical Engineering Department, University of California, Los Angeles, CA 90024, \EMAIL{ahujak@ucla.edu}} 
\AUTHOR{Simpson Zhang}
\AFF{Economics Department, University of California, Los Angeles, CA 90024, }
\AUTHOR{Mihaela van der Schaar}
\AFF{Electrical Engineering Department, University of California, Los Angeles, CA 90024, \EMAIL{mihaela@ee.ucla.edu}}
} 

\ABSTRACT{%
A newly launched online services firm typically desires to build a large user base for its website and must market itself to potential users that are not aware of its services. Once some users have entered the online firm, network effects, which are a result of features such as social interactions, content sharing etc., typically help to sustain the users' interest in the online firm. Marketing and network effects thus complement each other in the process of building the user base. This work is the first to study the interactions between marketing and network effects in determining the online firm's optimal marketing policy. We build a model in which the online firm starts with an initial user base and controls the growth of the user base by choosing the intensity of advertisements and referrals to potential users. A large user base provides more profits to the online firm, but building a large user base through advertisements and referrals is costly; therefore, the optimal policy must balance the marginal benefits of adding users against the marginal costs of sending advertisements and referrals. Our work offers three main insights: (1) The optimal policy prescribes that a new online firm should offer many advertisements and referrals initially, but then it should decrease advertisements and referrals over time. (2) If the network effects decrease, then the change in the optimal policy depends heavily on two factors i) the level of patience of the online firm, where patient online firms are oriented towards long term profits and impatient online firms are oriented towards short term profits and, ii) the size of the user base. If the online firm is very patient and if the network effects decrease, then the optimal policy prescribes it to be more aggressive in posting advertisements and referrals at low user base levels and less aggressive in posting advertisements and referrals at high user base levels. On the other hand, if the online firm is very impatient and if the network effects decrease, then the optimal policy prescribes it to be less aggressive in posting advertisements and referrals at low user base levels and more aggressive in posting advertisements and referrals at high user base levels. (3) The change in the optimal policy when network effects decrease also depends heavily on the heterogeneity in the user base, as measured in terms of the revenue generated by each user. An online firm that generates most of its revenue from a core group of users should be more aggressive and protective of its user base than a firm that generates revenue uniformly from its users.

}%



\maketitle

%


\section{Introduction}




%
%

%
\subsection{Motivation}

When an online services firm is launched it desires to build a large user base, but many potential users in the market may not be aware of the firm's website and services. Therefore, the online firm needs to reach out and market itself to potential users by posting advertisements and by giving referrals. As users become aware of the online firm they will enter the firm's website and then generate benefits for the firm by paying subscription fees, clicking on advertisements, etc \footnote{See http://finance.yahoo.com/news/must-know-assessingfacebook-
	revenue-170009607.html}. Users generate these benefits through the duration of their stay on the firm's website and hence the online firm desires not only to have a large user base but also to have its users remain with it for sufficiently long. Network effects, which result from social features like chatting with friends, user feedback etc., help to sustain each user's interest in the online firm and thus impact the length of time users continue to utilize its services. Network effects thus complement the force of marketing in the process of building a user base. However, no systematic models currently exist that can inform our understanding of how these two forces interact. For instance, the impact of the strength of network effects on the firm's optimal policy is not well understood, for instance, how should an online firm change its marketing decisions under different magnitudes of network effects. In this paper, we develop a general model that considers both marketing and network effects, and we determine the optimal marketing policy a firm should adopt to maximize its revenue.

\subsection{Model overview}
We develop a dynamic continuous-time model for the user base growth of an online services firm. The online firm starts with a small initial user base and markets to potential users by posting advertisements and by giving referrals. The potential users can also become aware of the online firm through other methods, for e.g. word-of-mouth, which are not controlled by the online firm. The online firm desires a large user base in order to procure more revenue, but building its user base through marketing is costly; therefore, it must design a policy that balances the marginal benefits of adding users against the marginal costs of marketing. We assume that the users arrive in the decreasing order of preference, which is proportional to the benefit a user generates for the website. This leads to online firm's revenue per unit time to be increasing and concave in the size of its user base, with the rate at which the benefit saturates reflecting the degree of heterogeneity in the user base. Users stay on the online firm for random durations, which due to network effects will increase in expectation as the user base size grows. The strength of the network effects is represented by the extent to which a growth in the user base influences the expected stay duration of each user. We assume that if the user base of an online firm increases then the number of friends of each user on the online firm should increase as well. The average time that a user interacts with each friend affects the increase in the stay duration of the user for every friend that is added and thus quantifies the strength of network effects in our model. Therefore, the user base dynamics of the users on the online firm are determined by three distinct factors: the intensity of direct arrivals through  word-of-mouth, the intensity of sponsored arrivals through advertisements, and the expected duration that each user stays on the online firm, which is impacted by the network effects. The online firm designs its policy taking future benefits into consideration, which depend on the user base dynamics of the online firm. The online firm discounts future benefits, with the extent of discounting depending on the duration the online firm expects to be in the market. 

\subsection{Results}
We show that if the online firm starts with a small initial user base, then the optimal policy is unique and prescribes the online firm to  give many advertisements and referrals initially and then decrease them over time. The optimal policy and the corresponding user base trajectory eventually converge to a unique steady state. Note that the optimal policy prescribes the online firm to decrease advertisements and referrals despite the fact that advertisements and referrals given at a later time can potentially be more beneficial; at later times, the user base is larger and each incoming user stays for a longer duration due to network effects. The factors that contribute to this result are: i) the marginal benefits decrease as the number of users increases, which is due to the concavity of the benefit function, ii) the stay time of the users is bounded above and saturates as the user base size increases.

Given the volatility of online settings, the strengh of network effects can change rapidly as new competitors or avenues for social interactions emerge, and so it is important to understand how the optimal policy is affected by a decline in the network effects. We show that if there is a decline in network effects, then the prescriptions of the optimal policy depend heavily on two factors: i) the online firm's patience level and ii) the size of the user base relative to an upper bound, which is an optimistic estimate of the online firm's user base in steady state. If there is a decline in network effects and if the online firm is very patient, then the optimal policy  will prescribe to be more aggressive in advertising  and sending referrals when the user base is very small relative to the upper bound, while it will prescribe  to be less aggressive in advertising  and sending referrals when the user base is close to the upper bound.  This result is reminiscent of the concept of barriers to entry in marketing. A very patient online firm operating at a small user base, which means that the market of the online firm is not saturated, increases the advertisements and referrals to create barriers to entry for the new online firms that cause a decline in its network effects. For instance, Netflix, which represents a very patient firm operating in a market (global market) that is not saturated, increased investment to expand its  user base \footnote{http://www.forbes.com/sites/petercohan/2015/01/21/4-reasons-to-invest-in-netflix/} and counter the increasing number of other rival firms.   
On the other hand, if there is  a decline in network effects and if the online firm is very impatient, then the  optimal policy will prescribe to be less aggressive in advertising and sending referrals when the user base is very small relative to the upper bound, while it will prescribe to be more aggressive in advertising and sending referrals when the user base is close to the upper bound.

The above comparison shows that the level of patience together with the size of the user base is crucial in determining the prescribed optimal policy. Interestingly, we also show that the level of heterogeneity in the user base measured in terms of the revenue generated by the users is also very important in determining optimal policy.
To the best of our knowledge, this is the first time that the impact of heterogeneity in the revenue generated by the user base on the optimal marketing policy has been systematically studied. We show that if there is a decrease in the network effects, the optimal policy prescribes an online firm with a very heterogeneous user base to increase its advertisements and referrals in the steady state. Examples of online firms that have a very heterogeneous user base are online gaming firms such as Zynga, where most of the users play for free and a very small fraction actually make in game purchases that contribute to revenue \footnote{ http://venturebeat.com/2014/02/26/only-0-15-of-mobile-gamers-account-for-50-percent-of-all-in-game-revenue-exclusive/}. Zynga  witnessed a decline in network effects \footnote{See alexa.com, https://en.wikipedia.org/wiki/Zynga} due to the rise of competing website King \footnote{http://www.forbes.com/sites/tomiogeron/2013/03/26/how-king-com-zoomed-up-the-social-gaming-charts/} and  increased investments \footnote{http://www.bloomberg.com/news/articles/2012-09-07/zynga-adds-more-developers-to-publishing-network-seeking-revenue} in response to the decline in network effects, which is similar to the prescriptions of the optimal policy in our model. On the other hand, we show that if there is a decrease in network effects, then the optimal policy prescribes a firm with a very homogeneous user base to decrease  its advertisements and referrals in the steady state. Examples of online firms that have a very homogeneous user base are subscription based online firms such as Netflix, where all the users need to pay the subscription fee to use the service.


In practice the online firm cannot perfectly predict the user arrivals and exits since the users' behavior is uncertain. We analyze how should the prescriptions of the optimal policy change to account for this uncertainty. We model the uncertainty in the user arrivals per unit time using a standard geometric Brownian motion process, which captures the uncertainty of user arrivals and  exits.  We show that if there is high uncertainty in user arrivals, then the optimal policy will prescribe to be less aggressive in advertising and sending referrals.  We also show that if there is higher uncertainty in the user arrivals, then the optimal policy will achieve a smaller expected user base in steady state as well as achieve lower expected profits  compared to the deterministic case with no randomness in arrivals/exits. Basically, this happens because the increase in the random arrivals cannot compensate for the increase in the random exits and the decrease in the sponsored arrivals.

\subsection{Key findings}

We summarize the key findings of this paper below.
\begin{enumerate}
	\item  If the online firm starts with a small initial user base, then  the optimal policy is unique and it prescribes the online firm  to be more aggressive in advertising and send more referrals initially and then decrease the advertisements and referrals over time.  
	\item If there is a decline in network effects, then the  prescriptions of the optimal policy  depend jointly on the level of patience and the  size of the user base. The optimal policy prescribes a  very patient (impatient) online firm  to increase (decrease) the advertisements and referrals at small user base levels and decrease (increase) the advertisements and referrals at large user base levels. 
	\item If there is a decline in the strength of the network effects, then we also show that the prescriptions of the optimal policy  depend on the level of heterogeneity in the user base. The optimal policy prescribes an online firm with very heterogeneous (homogeneous) user base  to increase (decrease) the advertisements and referrals in steady state.
	\item If the user arrivals/exits are uncertain, then we show that the optimal policy prescribes the online firm to be less aggressive in comparison to the case when the user arrivals/exits are deterministic. 
\end{enumerate}

The rest of the paper is organized as follows. In Section 2 we describe the related literature. In Section 3 we present the model and the policy design problem of the online firm. In Section 4 we analyze the existence and uniqueness of the optimal policy. In Section 5 we analyze the behavior of the optimal policy on the path to steady state. In Section 6 we analyze the impact of change in network effects on the optimal policy. In Section 7 we discuss the extensions of the model and in Section 8 we conclude the paper.

\section{Related Literature}

\subsection{Literature on advertising}
In this section we summarize the relevant works on advertising. The main difference between  these works and our work is that these works do not incorporate network effects into their model and thus do not analyze the impact of network effects on advertising.

The seminal work of  \cite{ozga1960imperfect} was the first to study the lack of perfect knowledge about the market among buyers and sellers. It developed an analytical model of advertising based on information diffusion amongst the buyers. \cite{nerlove1962optimal}   considered a more general setting by incorporating the impact of advertising expenditures on the demand of the product into a firm's decision making and proposed optimal advertising and pricing policies. The 
Bass model (\citealt{doi:10.1287/mnsc.15.5.215}) was the first to analyze the process of the adoption of goods in a market at a micro level. It  separated users in the market into two types, innovators and imitators and the two types followed different product adoption dynamics. The Bass model did not separately incorporate decision making such as the intensity of advertising and pricing in the market; however, the model has been tremendously successful in capturing the product adoption dynamics in several scenarios as discussed in \cite{Bass2}. 
 \cite{robinson1975dynamic} separately incorporated pricing decisions  of the firm into Bass's model. They considered patient firms that derive pricing policies taking the evolution of the market, i.e. demand,  as predicted by Bass's model into consideration.  \cite{thompson1984optimal} then considered a more general model with patient firms that derive joint pricing and advertising policies taking the evolution of the market into consideration.

The above models are well suited for a market where the users generate one time benefits for the firm upon purchase of the product and leave. However, the rise of service-based economies has shifted the focus to scenarios where the relationship of the firm with the user is long-term in nature (for e.g., online firms, banks etc.) and this has been a subject of investigation recently  (see for instance \citealt{gupta2006modeling},  \citealt{gourio2014customer}).  \cite{gupta2006modeling} discusses various  metrics based on different models that are used to evaluate the value of a user's stay-time with the firm. The work compares these different metrics based on empirical insights.  

 Theoretical research analyzing the dynamics of a  user base and its impact on the firm's policies has been limited.  \cite{gourio2014customer}  analyzes the impact of search frictions in a market arising due to inefficiencies in acquiring users for the firm. They cite the enormous spending done by the firms for marketing and selling as  evidence of such inefficiencies. They  provide a micro-foundation for the search process and analyze the impact of search frictions on a firm's profits, values, sales and markups. Our work makes several similar modeling assumptions i) the user's relationship with the online firm is long-term, ii) the online firm's benefits depend on the stay duration of the users, and iii) the process of advertising is inefficient, which causes the cost of advertising to be convex. A key  difference between our model and all of the above works is that we consider the advertising policy of the online firm in conjunction with network effects, i.e. larger networks positively impact the user's stay duration on the firm's website. Moreover, the questions that we investigate are different from the existing works and concern the impact of different network effect strengths on the optimal policy. Next, we discuss the relevant works which study the impact of different network effects on firm's policies, especially with respect to pricing.

\subsection{Literature on network effects} 
In this section we summarize the relevant works that study network effects. The main difference between  these works and our work is that these works do not incorporate advertising into their model and thus do not analyze the impact of network effects on advertising.

The seminal work of \cite{katz} was the first to formally model and analyze the issue of pricing in various market scenarios where network effects play an important role. Since then, the issue of pricing under network effects has been addressed in different works analyzing several different scenarios (see for instance \citealt{candogan1}, \citealt{lobel}). Most of the works (e.g., \citealt{katz}, \citealt{candogan1}, \citealt{lobel})  consider a static setting but
some of the recent works (e.g., \citealt{markovich}, \citealt{cabral}) consider a dynamic setting where the users choose which network to enter while taking the price set by the network, the preference for the network and the size of the network into consideration.  \cite{fudenberg1999pricing} consider the dynamics of an incumbent and an entrant  in a market with network externalities, and they show that the incumbent can set lower prices  in the equilibrium due to the  competition faced from the entrant. In each of these works, the authors do not model the advertising policy of the online firm, which will impact the rate at which users become aware of the network. To the best of our knowledge there is a gap between the literature on advertising and network effects. In this work we take an  important first step to bridge this gap.

%
%
%

%
\section{Model}

In this section we propose a model to understand how aggressively an online services firm should post advertisements and  send referrals in order to influence its user base dynamics. The two key features of the model are: i) the user base
dynamics are impacted by network effects and ii) the users generate different revenues for the online firms, which causes the user base to be heterogeneous.

\subsection{User arrivals and exits}

We consider a setting with a continuum of potential users and an online firm
which markets itself to these users by posting advertisements on other online firms (e.g., Facebook, Google etc.) or by giving referrals. New
users visit the online firm either through sponsored media (advertisements/referrals) or arrive directly
through exogenous methods, such as directly visiting the online firm or arriving at the online firm after hearing about it from others. We refer to these exogenous arrivals as direct arrivals and denote them as $f(x)$, where $x$ is the size of the user base of the online firm. Note that they are not directly influenced by the online firm. We assume that  $f(x)$ is increasing, non-negative, continuously differentiable and  bounded above and below by positive constants
$f^{sup}$ and $f^{inf}$, respectively. $f(x)$ is increasing because when the user base is large  more users become aware of the online firm and thus enter it. 

Each user who visits the online firm will stay on it and
then eventually exit the online firm after some random time, which is determined by a Poisson arrival process that starts along with the arrival of the user. At the first arrival of the process the user leaves. The exit can be interpreted, for instance, as the user losing interest in the online firm. 
The rate of Poisson arrivals is given by $\frac{1}{\tilde{\eta}+\eta g(x)}$ and hence, the stay
time has an exponential distribution with a mean given by $\tilde{\eta}+\eta g(x)$. The first term
$\tilde{\eta}$ represents the expected stay time in the absence of the
network effects; it is the average time a user will spend
independent of the presence of others. The second term $\eta g(x)$ represents the expected stay time due to the network effects; it is the average time a user will spend on social interactions. The term $\eta$ represents the average time that a user spends per interaction and hence, it is a measure of the strength of network effects. Note that $\eta$ is a fixed constant that depends on the type of the online firm (for e.g., Netflix features much less social interactions than Facebook) and the presence of other available venues for social interaction.  We assume that $\eta$ is bounded, i.e. $0 \leq \eta\leq \eta^{sup}$. The term $g(x)$ represents the average number of interactions that a user carries out. We assume that $g(x)$ is non-negative, increasing, continuously differentiable and bounded above by $g^{sup}$. 
$g(x)$ is increasing  because when the user base is large each user will have more opportunities to interact.  We denote the size of the user base that is using online firm's services at time
$t$ as $x(t)$, where $x:[0,\infty)\rightarrow[0,\infty)$ is the user base trajectory of the users.  We restrict our attention to user base trajectories that are absolutely continuous. The initial
user base $x(0)$ is denoted by $x(0)=x_{initial}$, which is assumed to be  small, $0<x_{initial} <  f^{inf}\tilde{\eta}$.


\subsection{Online firm's benefits and costs}

The size of the user base $x$ is proportional to the unique user
statistic, which is often used as a metric to evaluate an online firm's
popularity and hence, it is critical to the revenue generated by the
online firm \footnote{http://www.pcmag.com/encyclopedia/term/53438/unique-visitors}. We assume that each potential user has a fixed preference for the online firm that is proportional to the benefit that the user will generate for the online firm.  We assume that the users enter the online firm in a decreasing order of their preference. We can show that the size of the user base $x$ on the online firm is sufficient to determine the preference distribution of the users on the online firm and thereby the benefit function (see  Appendix J).  Hence, the benefit generated by the online firm per unit time with a user base of size $x$ is denoted as $b(x)$, where $b:\mathbb{R}\rightarrow[0,\infty)$
 is a continuously differentiable increasing function. In fact it can also be shown that the benefit function $b(x)$ is concave (see  Appendix J).  The rate at which $b(x)$ saturates reflects the heterogeneity in the user base in terms of the revenue generated by the users. Next, we give an example of a concave benefit function and show that the rate of saturation of $b(x)$ reflects the extent of heterogeneity in the user base. 

\textit{Example-} Consider an online firm and a user base of potential users, which has a total mass of 1. Each user has a preference for the online firm $\omega$ drawn from a uniform distribution $[0,1]$ independent of others. Preference $\omega$ represents the willingness to pay for services of the firm. A user with preference $\omega$ will generate a benefit of $\frac{a}{(1-\omega)^{1-a}}$ for the online firm upon joining the online firm, where $0 \leq a\leq 1$.  If $a$ is close to $1$, then the user base is homogeneous, i.e. users generate similar benefits for the online firm. For example, subscription based firms such as Netflix have the same subscription fee for all the users. If $a$ is close to $0$, then the user base is heterogeneous, i.e. different users generate different benefits. For example,  online gaming firms such as Zynga have significant heterogeneity in benefits generated by users. Only 0.15 percent of the users generate 50 percent of the revenue for Zynga (See  \footnote{http://venturebeat.com/2014/02/26/only-0-15-of-mobile-gamers-account-for-50-percent-of-all-in-game-revenue-exclusive/}). We assume that the users join the online firm in the order of their preferences. We also assume that each user who joins the online firm will not leave. Note that this assumption is made for a better exposition and can be relaxed.  The total benefit  generated by a user base of $x$ users, where $x<1$ is $b(x) = \int_{1-x}^{1} \frac{a}{(1-\omega)^{1-a}}d\omega =-(1-\omega)^{a}]_{1-x}^{1}= x^{a} $. $\hfill \blacksquare$

 The long-term benefit of the online firm considering a discount
rate of $\rho$ can be computed as $B(\bar{x})=\int_{0}^{\infty}b(x(t))e^{-\rho t}dt$, where $\bar{x}$ is a succinct representation for the entire trajectory of $x(t)$ for all $t\geq 0$. Note that $\frac{1}{\rho}$ represents the discount factor and is proportional to the expected duration of stay of the online firm in the market. We refer to a  online firm with a high discount factor as patient and an online firm with a low discount factor   as  impatient.  

The online firm chooses the intensity at which it advertises itself and sends referrals and as a result it controls the
rate of sponsored arrivals $\lambda(t)$, where $\lambda:[0,\infty)\rightarrow[0,\lambda^{sup}]$
is an absolutely continuous and Lebesgue measurable function. The online firm bears
a cost per unit time for the advertisements, which increases with
the intensity of sponsored user arrivals $\lambda(t)$. It is given
by $c(\lambda(t))$, where $c:\mathbb{R}\rightarrow[0,\infty)$ is a
continuously differentiable and increasing function. The online firm needs
to post a large number of advertisements in order to find users interested in using its services.  Hence,
we assume (as in \citealt{gourio2014customer}) that the cost $c(\lambda)$
is strictly convex in $\lambda$.
The long-term
discounted average cost to the online firm is $C(\bar{\lambda})=\int_{0}^{\infty}c(\lambda(t))e^{-\rho t}dt$, where $\bar{\lambda}$  is a succinct representation for the entire policy $\lambda(t)$ for all $t\geq 0$.
The overall rate of change of the user base on the online firm is
\begin{eqnarray*}
\frac{dx}{dt}=f(x)+\lambda(t)-\frac{1}{(\tilde{\eta}+\eta g(x))}x
\end{eqnarray*}
The first two terms in the differential equation represent the rate
of direct and sponsored arrivals respectively and the third term represents the rate at which the users leave the online firm. Note that the third term is deterministic despite the fact that each user exits randomly. This is due to the fact that there is a continuum of users and each user will leave the online firm with probability $\frac{1}{\tilde{\eta} +\eta g(x)}\Delta t$ in the time interval $\Delta t$, with $\Delta t$  small, independent of other users. Hence, the number of users who leave the online firm in $\Delta t$ time is $\frac{1}{\tilde{\eta} +\eta g(x)} x\Delta t$.
 Table 1 
summarizes the notation  used in this paper and classifies
the variables based on the type of control that the online firm has on each
of them  (direct, indirect or no control).

While in this paper we assume that  the users arrive in the decreasing order of their preference, which results in a concave benefit function, we  also extend the model and  the results to incorporate a different order of arrivals of the users, in which the benefit functions are not necessarily concave  (see Appendix G). We can also relax the assumption that the users are characterized by a decision function that leads to the direct arrivals at the rate $f(x)$ and sponsored arrivals at the rate $\lambda(t)$. We extend the model and results to scenarios where the users take decisions to enter the online firm or not based on a utility function that depends on the preference of the user and the size of the user base (see  Appendix H). 
 
In the next section we describe the policy design problem of the online firm.

\begin{table}
	\centering
	\begin{tabular}{ | l | l | l| }
		\hline
		Variables & Notation   & Control type \\ \hline
		Expected stay time at time t & $\tilde{\eta}+\eta g(x(t))$  & Indirect \\ \hline
		User base size at time t &$x(t)$ & Indirect \\ \hline
		Direct arrivals at time t & $f(x(t))$  & Indirect\\ \hline
		Expected stay time (network effects part) at time $t$ &  $\eta g(x(t))$ & Indirect \\ \hline
		Average number of social interactions per user & $g(x(t))$ & Indirect \\ \hline 
		Flow benefit at time $t$  & $b(x(t))$ &  Indirect  \\ \hline 
		Long-term discounted average benefit & $B(\bar{x})$  &  Indirect\\  \hline
		Sponsored arrivals at time t & $\lambda(t)$   & Direct \\ \hline
		
		Flow costs at time $t$ & $c(\lambda(t))$ & Direct \\ \hline
	
		Long-term discounted average cost   & $C(\bar{\lambda})$  & Direct \\ \hline
		Network effect strength&  $\eta$  & No control \\ \hline 
		Expected stay time (independent part) &  $\tilde{\eta}$ & No control \\ \hline
		Initial user base & $x_{initial}$ & No control \\ \hline
		Discount rate   & $\rho$ & No control \\ \hline
		Upper bound on direct arrivals & $f^{sup}$& No control \\ \hline
		Lower bound on direct arrivals & $f^{inf}$ & No control \\ \hline
		Upper bound on strength of network effects & $\eta^{sup}$ & No control \\ \hline
		Upper bound on social interactions per user & $g^{sup}$ & No control \\ 
		\hline 
	\end{tabular}
\caption{Classifying the variables in the model based on the control type}
\end{table}


\subsection{Online firm's policy design problem}

 The online firm desires to maximize its long-term discounted profit, i.e. $B(\bar{x})-C(\bar{\lambda})$.  It solves for the optimal policy that maximizes its long-term discounted profit, which is formally stated below.  Recall that the initial user base $x(0)=x_{initial}$ is given

 \begin{eqnarray*}
&& \textbf{Policy Design Problem}\\
 &&  \max_{\bar{\lambda}} B(\bar{x})-C(\bar{\lambda})= \int_{0}^{\infty}(b(x(t)) -c(\lambda(t)))e^{-\rho t}dt \\
 && \text{subject to}\; \frac{dx}{dt} = f(x) + \lambda(t) - \frac{1}{(\tilde{\eta} +\eta g(x))}x\\
 \end{eqnarray*}
 We denote a policy that solves the above continuous-time optimization problem  $ \lambda_{x_{initial}}^{*}$, where $\lambda_{x_{initial}}^{*}:[0, \infty) \rightarrow [0,\lambda^{sup}]$ and we refer to it as the optimal policy as a function of time. The corresponding optimal user base trajectory is denoted as $x_{x_{initial}}^{*}$, where $x_{x_{initial}}^{*}:[0, \infty) \rightarrow [0,\infty]$ is a function which maps each time instance $t$ to the user base at that time. 

The maximum value of the long-term profit that is achieved by solving the above  continuous-time optimization problem is denoted as $\Pi^{*}(x_{initial},\eta)$. $\Pi^{*}(x_{initial},\eta)$  has to satisfy the Hamilton-Jacobi-Bellman (HJB) equation, which is given below, provided $\Pi^{*}(x_{initial},\eta)$  is differentiable with respect to $x_{initial}$ (see  \citealt{kushner1971introduction}). Next, using the HJB equation we will compute an optimal policy that will  depend only on the size of the user base and  relate it with the optimal policy as a function of time $ \lambda_{x_{initial}}^{*}$. This approach is standard in optimal control (See \cite{speyer2008stochastic}).

\begin{eqnarray*}
\rho \Pi(x,\eta) =  \max_{0 \leq \lambda \leq \lambda^{sup}}\big(b(x)) -c(\lambda) + \frac{\partial \Pi(x,\eta)}{\partial x}(f(x)+\lambda- \frac{1}{\tilde{\eta}+ \eta g(x)}x)\big)
\end{eqnarray*}
We substitute the  optimal solution  of the policy design problem $\Pi^{*}(x_{initial},\eta)$  in the HJB equation. The maximizer of the right hand side in the HJB equation is given as $\zeta^{*}(x_{initial}) =\arg \max_{0 \leq \lambda \leq \lambda^{sup}}\big(b(x_{initial})) -c(\lambda) + \frac{d\Pi^{*}(x_{initial},\eta)}{dx}(f(x_{initial})+\lambda- \frac{1}{\tilde{\eta}+ \eta g(x_{initial})}x_{initial})\big)$, where $\zeta^{*}(x_{initial})$  corresponds to the optimal policy  at the initial user base size $x_{initial}$. This follows from the Bellman's principle of optimality (see \citealt{kushner1971introduction}). Note that this optimal policy $\zeta^{*}(x_{initial})$ only depends on the initial size of the user base. Following the same method we can obtain an optimal policy at different user base sizes $x$ using the optimal solution $\Pi^{*}(x,\eta)$. We can thus define an optimal policy as a function of the user base size, denoted as $\zeta^{*}$, where $\zeta^{*}:[0,\infty) \rightarrow [0,\lambda^{sup}]$ is a mapping from the user base size $x$ to the optimal policy. If  the optimal policy as a function of the  user base size $\zeta^{*}$ and the  optimal policy as a function of time  $\lambda^{*}_{x_{initial}}$ are unique, then $ \lambda^{*}_{x_{initial}}(t) = \zeta^{*}(x^{*}_{x_{initial}}(t))$  has to hold for consistency.  

   In the next section we will discuss the existence and uniqueness of  the optimal policy  as a function of the user base and the optimal policy as a function of time.

\section{Existence and uniqueness of the optimal policy}
In this section we first show that the optimal policy as a function of user base $\zeta^{*}$ exists and is unique. Next, we show that the optimal policy  as a function of time $\lambda_{x_{initial}}^{*}$ also exists and is unique. 

First, we state some assumptions that are required in the rest of the paper.  Note that it can be shown that  every user base trajectory satisfying the differential equation in the policy design problem is always bounded above by $x_{u}$ (see  Appendix A for the expression)  and bounded below by $x_{initial}$  (see  Appendix A for the proof) \footnote{The population trajectory not falling below $x_{initial}$ is a result of the deterministic model in Section 3. In Section 7 we work with a stochastic model, where the population trajectory can fall below the initial population level.}.   Therefore, in the next few assumptions we only require the properties stated to hold in an interval $[x_{initial},x_{u}]$.
\begin{assumption}
	\begin{enumerate}
\item	$c^{'}(0)<\frac{\partial \Pi^{*}(x,\eta)}{\partial x}<c^{'}(\lambda^{sup}),\;\forall x \in [x_{initial},x_{u}]$, where $c^{'}(\lambda) = \frac{dc(\lambda)}{d \lambda}$ and $\frac{\partial \Pi^{*}(x,\eta)}{\partial x}$ is the partial derivative of $\Pi^{*}(x,\eta)$ with respect to $x$. 
	
\item	$\frac{\partial^{2} \Pi^{*}(x,\eta)}{\partial x^{2}}$ exists and is continuous with respect to $x$, $\forall x\in [x_{initial}, x_{u}]$.
	\end{enumerate}	
\end{assumption}

We can show that $\frac{\partial \Pi^{*}(x,\eta)}{\partial x}$ exists (see  Appendix A).
The first part of Assumption 1 requires  $\frac{\partial \Pi^{*}(x,\eta)}{\partial x}$ to be bounded above and below,  which should hold for  online firms  whose profits  do not change  too sharply with a change in the size of the user base $x$. This assumption ensures that the optimal policy is always in the interior of $[0,\lambda^{sup}]$. 
 The second part of the assumption requires that $\frac{\partial^{2} \Pi^{*}(x,\eta)}{\partial x^{2}}$ exists and is continuous, which should hold for  online firms whose marginal increase in profits as users are added does not change too sharply with an increase in the size of the user base $x$.   We discuss in Appendix A cases when   Assumption 1 does not hold. In these cases the optimal policy will not necessarily be in the interior of the set $[0,\lambda^{sup}]$, despite that we can show that the flavor of the main results in this paper will not change.
%

We define the rate of arrivals with no advertisements and referrals as follows
$h_{\eta}(x)=f(x)-\frac{1}{(\tilde{\eta}+\eta g(x))}x$ \begin{assumption}  $h_{\eta}(x)$ is  continuously differentiable and concave in the interval $[x_{initial},x_{u}]$.
	
%
	
	 \end{assumption}

%

In Assumption 2 we assume that the rate of arrivals in the absence of advertisements and referrals $h_{\eta}(x)$ is concave. 
This implies that the marginal increase in the arrivals with an increase in the user base decreases at larger user base levels. This assumption holds for online firms
for which the popularity increases faster at relatively small user base
levels and gradually saturates as the online firm becomes larger. 

Having stated the assumptions we now move to the main results of this section.

\begin{proposition}
	The optimal policy as a function of user base $\zeta^{*}$ exists and is unique.
	\end{proposition}

For all the detailed proofs refer to Appendix given at the end. The proof of the above proposition is given in Appendix A.


 \begin{proposition} The  optimal policy as a function of time $\lambda_{x_{initial}}^{*}$ and
	the corresponding user base trajectory as a function of time $x_{x_{initial}}^{*}$ exists and is unique. \end{proposition}

The proof of the above proposition is given in Appendix A.

Note that Proposition 2 does not follow directly from Proposition  1. To show Proposition 2 we need to prove that once $\zeta^{*}$ is uniquely determined  there exists a unique solution to $\frac{dx}{dt} = f(x) + \lambda(t) - \frac{1}{(\tilde{\eta} +\eta g(x))}x$ for a given initial condition $x(0)=x_{initial}$, which will be denoted as $x^{*}_{x_{initial}}$.

 Now we know that both the policy as a function of the user base $\zeta^{*}$ and the policy as a function of time $\lambda_{x_{initial}}^{*}$ exist and are unique.  The uniqueness of the policy is important in order to compare the change in the optimal policy under varying network effects. Note that as we already pointed out the following relation will hold  $\lambda_{x_{initial}}^{*}(t) = \zeta^{*}(x^{*}_{x_{initial}}(t))$. Also, when the benefit and cost function are quadratic and $f(x) + \lambda(t) - \frac{1}{(\tilde{\eta} +\eta g(x))}x$ is linear in $x$, then a closed form expression can be obtained for the optimal policy as a function of time and optimal policy as a function of user base size (See Appendix I). 
 
 In the next section we discuss the convergence  of the optimal policy 
 and also investigate the questions concerning how  the optimal policy changes as a function of time.


\section{Optimal policy behavior and convergence to steady state}

In this section we  analyze the convergence of the optimal policy as a function of time $\lambda_{x_{initial}}^{*}$ and the corresponding optimal user base trajectory $x_{x_{initial}}^{*}$. In order for us to analyze convergence
we define the steady state of the policy as follows. \begin{definition}
	The steady state of a system is a pair $(\lambda_{s},x_{s})$, which satisfies the following condition: if the initial user base in the policy design problem is $x_{s}$,
	then $\lambda(t)=\lambda_{s},\;\forall t$ is the optimal policy as a function of time
	and the corresponding user base trajectory is $x(t)=x_{s}\;\forall t$.
\end{definition}

The above definition of the steady state states that if the
initial user base is at the steady state value then the
optimal policy will correspond to the exact number of sponsored arrivals such
that the user base stays at that fixed value. The next theorem proves that the steady state exists and is unique.

%
%
%

\begin{theorem}
	The steady state $(\lambda_{s}, x_{s})$ of the online firm exists and is unique. 
\end{theorem}
For the proof of the above theorem refer to Appendix B given at the end. 

Theorem 1 shows that the steady state exists and is unique. If we can show that the optimal policy converges to the steady state, then understanding the impact of changing network effects on the steady state,  which is unique, will serve as an important tool for analysis in the next sections. 
%
Next we prove that the optimal policy converges to the steady state.
%


\begin{theorem} The optimal policy as a function of time $\lambda_{x_{initial}}^{*}(t)$ and the corresponding optimal user base trajectory as a function of time $ x_{x_{initial}}^{*}(t)$ converge to the steady state values, i.e. $\lim_{t\rightarrow \infty}\lambda_{x_{initial}}^{*}(t) \rightarrow \lambda_{s}$ and $\lim_{t\rightarrow \infty} x_{x_{initial}}^{*}(t) \rightarrow  x_{s}$.
\end{theorem}

For the proof of the above theorem refer to Appendix C given at the end. 

The intuition behind the proof of the above theorem is straightforward.  We know that the  initial user base is smaller than the smallest steady state user base level $f^{inf}\tilde{\eta}$. Based on this it can be argued that the optimal user base trajectory is always increasing. We also know that all the user base trajectories are bounded above hence,  the optimal user base trajectory will converge. Next, we analyze the behavior of the optimal policy on its path to the steady state: When the online firm is launched, should it give more referrals in the beginning? Or should it give more referrals at a later stage?. The online firm faces the following trade-off in choosing among these two approaches. The first approach requires the online firm to market aggressively in the beginning and achieve a large user base quickly. A large user base will generate more benefits for longer durations but requires the online firm to pay more costs for marketing. The second approach requires the online firm to wait for the user base to increase sufficiently  and then market aggressively. The second approach allows the online firm to limit costs and invest later when the returns from each advertisement and  referral  will be high since the users will stay longer due to the network effects. The next theorem determines which approach is better.

\begin{theorem} The optimal policy as a function of time $\lambda_{x_{initial}}^{*}(t)$  decreases with time $t$ and the corresponding user base trajectory as a function of time $x_{x_{initial}}^{*}(t)$ increases with time $t$.
\end{theorem}

For the proof of the above theorem refer to Appendix D.

We show that the marginal increase in  profit from adding users at a larger user base  is lower and thus conclude that the first approach will be superior to the second. In the second approach, the online firm delays the choice of giving more referrals to achieve larger gains at a later stage. Due to the decrease in the marginal increase in  profit from adding users, the online firm cannot  compensate for the loss it has to bear because of a smaller user base in the earlier stages.

In the next section we discuss the impact of the network effects on the optimal policy.

\section{Impact of a change in network effects on the optimal policy}

 How does a decrease in the strength of network effects impact the prescriptions of the optimal policy to the online firm? A decrease in the strength of network effects can result from an increase in other avenues for social interactions   or the emergence of other online services  which leads to an increase in indirect competition. If there is a decrease in network effects, then the following two comparisons are relevant:
 \begin{enumerate}
 	\item Comparing  the prescriptions of the optimal policy $\zeta^{*}(x)$ at the same user base size $x$. 
 	\item Comparing  the prescriptions of the optimal policy $\lambda^{*}_{x_{initial}}(t)$ at the same time $t$. 
 \end{enumerate}
 
  A comparison at the same user base size is significant because the same user base size corresponds to same level of popularity for the online firm, while  a comparison at the same time is significant because  the same time corresponds to the same date after the launch of the online firm. 


\subsection{Impact of the network effects: comparing the optimal policy at the same size of the user base}

In this section we analyze the impact of a decrease in the strength of network effects on the optimal policy $\zeta^{*}(x)$ at the same size of the user base. Next, we state some assumptions needed exclusively for the results in this subsection (Theorem 4 and Corollary 1). 

 We assume that $c^{'}(\lambda)$ has a positive lower bound $c^{inf}$. $c^{'}(\lambda)$  measures the additional cost that the online firm needs to bear for increasing the sponsored arrival rate by a unit value.    The online firm will post more advertisements to increase the sponsored arrivals. Typically, there is a minimum cost for posting an advertisement, which justifies this assumption. We also assume that $\Pi^{*}(x,\eta)$ is continuously differentiable with respect to $\eta$ in the interior of the interval $[0, \eta^{sup}]$.  In addition, we also assume that the rate of arrivals at $x_{initial}$ in the absence of referrals and advertisements is sufficiently high $h_{\eta}(x_{initial})> \frac{dv(\lambda^{sup})}{d\lambda}, \; \forall \eta \in [0,\eta^{sup}]$, where $v(\lambda)=c(c'^{-1}(\lambda))$ and $c'^{-1}(\lambda)$ is the inverse of the function $c^{'}(\lambda)$.

 For the next theorem we define a lower bound $\rho^{l}$ and an upper bound $\rho^{u}$ on the discount rate, where $\rho^{l}<\rho^{u}$. If $\rho<\rho^{l}$ then the online firm is very patient and if $\rho > \rho^{u}$ then the online firm is very impatient. We also define a lower bound $x_{1}$ and an upper bound $x_{2}$ on the user base size, where $x_{initial} \leq x_{1} \leq x_{2} \leq x_{u}$.
 The expressions for these thresholds are derived and  stated in Appendix E.

\begin{theorem}
	If there is  a decrease in the strength of network effects $\eta$ and
	\begin{itemize}
		\item  if the online firm is very patient ($\rho \leq \rho^{l}$) and 
	 has a sufficiently small user base $ x\leq x_{1}$,  then the optimal policy prescribes  to increase the advertisements and referrals;
		\item  if the online firm is very patient ($\rho \leq \rho^{l}$) and has a sufficiently large user base $x \geq x_{2}$, then the optimal policy prescribes  to decrease the advertisements and referrals;
			\item  if the online firm is very impatient ($\rho \geq \rho^{u}$) and 
			has a sufficiently small user base $ x\leq x_{1}$,  then the optimal policy prescribes to decrease the advertisements and referrals;
			\item  if the online firm is very impatient ($\rho \geq \rho^{u}$) and has a sufficiently large user base $x \geq x_{2}$, then the optimal policy prescribes to increase the advertisements and referrals.
	\end{itemize}
\end{theorem}

For the proof of the above theorem refer to Appendix E. 

Before we discuss the above theorem, we point out that the result of the above theorem does not depend on the assumption that the benefit function is concave and this is discussed in detail in Appendix G.

In the above theorem we show that the change in the optimal policy in response to a change in the network effects depends on the size of the user base and the level of patience. We first discuss the intuition for the result on very patient online firms. Examples of very patient online firms include subscription based online firms as these online firms  provide services and expect the users to continually use them.  If there is a decrease in  network effects, then it can potentially lead to a reduction in the user base size in the future, causing an increase in the marginal benefits from adding users. If the user base is sufficiently small then the increase in the marginal benefits is large and since the online firm is very patient the optimal policy prescribes giving out more referrals and advertisements such that it can achieve a large user base in the long run that compensates for the costs  of giving more referrals and advertisements.  On the other hand, if there is a decrease in network effects and if the online firm has a sufficiently large user base, which is close to its steady state level, then it does not expect to achieve a much larger user base than its existing user base in the long run. Therefore, the optimal policy  prescribes decreasing the advertisements and referrals and saving the costs instead. Next, we discuss the intuition for the   result with very impatient online firms.

Examples of very impatient online firms are  online gaming firms because  games typically have a short life-time. 
If there is a decrease in  network effects and if the online firm has a small user base, then potentially there will be a reduction in the user base size in the future.  The online firm is impatient and will not increase the advertisements and referrals because increasing advertisements and referrals will lead to an increase in the user base only in the long-term, while the online firm values short-term benefits. Also, an increase in the  advertisements and referrals will cause a large increase in the short-term costs. This is because when the user base is small the sponsored arrival rate is already high and increasing it further will be very costly. On the other hand, the opposite would be true if the online firm has a sufficiently large user base. At large user base  the rate of sponsored arrivals is low hence, the marginal costs of advertising are low as well.  Thus, the online firm increases  advertisements and referrals to maintain the  large user base that generates sufficient benefits in the short-run and compensates the costs of the increase in referrals and advertisements. 
The results of the above theorem are summarized in Table II. 
%

\begin{table}
	\centering
	\begin{tabular}{ |l| l | l| }
		\hline
		Firm type based on $\rho$ & Size of the user base  & Change in referrals and ads \\
		\hline
		Very patient & Small & Increase \\ 
		\hline
		Very patient & Large & Decrease \\
		\hline 
	Very impatient & Small & Decrease \\
		\hline 
Very impatient & Large & Increase \\
		\hline
	\end{tabular}
\caption{Impact of the decrease in the network effects on the referrals and advertisements}
\end{table}

In this section we learned how  online firms should change their optimal policy in reaction to a reduction in the network effects. We saw that the the prescriptions of the optimal policy crucially depend on the online firm's level of patience and the size of the user base.
Next, we state a corollary of Theorem 4 that compares  the optimal policy $\lambda^{*}_{x_{x_{initial}}}$ at time $t=0$ under the impact of decreasing network effects. For this corollary we assume that $x_{initial}<x_{1}$ in addition to the assumption in Section 3.1 that $x_{initial}<f^{inf}\tilde{\eta}$.

   \begin{corollary}
   	 If there is a decrease in the strength  of network effects $\eta$ and 
   	\begin{itemize}
   		\item if the online firm is very patient ($\rho \leq \rho^{l}$), then the optimal policy prescribes  to increase the advertisements and referrals at initial time instances close to $t=0$;
   		\item if the online firm is very impatient ($\rho \geq \rho^{u}$), then the optimal policy prescribes to  decrease the advertisements and referrals at initial time instances close to $t=0$.
   	\end{itemize} 
   \end{corollary}
   
   At time $t=0$ the population trajectory will have the same value $x_{initial}$ and therefore, we can use the comparison at low population levels in Theorem 4 to derive the above Corollary.  In the next section we compare the optimal policy at large times $t\rightarrow \infty$ under the impact of decreasing network effects.

\subsection{Impact of the network effects: comparing the optimal policy at the same time}


In the previous section we already compared the change in the prescriptions of the optimal policy at time $t=0$ (see Corollary 1). Next, we will analyze the impact of a change in network effects on the optimal policy at time $t\rightarrow \infty$, i.e. in the steady state. The results that we will discuss extend to sufficiently large values of time instances (due to the continuity properties of the optimal policy).  Recall that in Section 2 we proved that  the optimal policy will converge to the steady state  denoted as $(\lambda_{s},x_{s})$. In essence, we will analyze how  $\lambda_{s}$ changes when the strength of network effects decreases. 

\begin{theorem}
	If there is a decrease in the strength of network effects $\eta$ and 
	\begin{itemize}
		\item if the online firm is very patient ($0 \leq \rho \leq \rho^{l}$), then the optimal policy prescribes  to  decrease the advertisements and referrals in the steady state;
		\item if the online firm is very impatient ($\rho \geq \rho^{u}$), then the optimal policy prescribes to  increase the advertisements and referrals in the steady state.
	\end{itemize}
	\end{theorem}
	
	For the proof of the above theorem refer to Appendix F. 
	
 If the online firm is very patient, then a decrease in the network effects causes the online firm to increase the advertisements and referrals in the initial stages close to $t=0$ (see Corollary 1) such that it can build a large user base  and derive more benefits in the long run. As  time goes on and the user base becomes sufficiently large the online firm decreases the advertisements and referrals and saves the costs as it does not expect the user base size to increase much.  
%
 On the other hand if the online firm is very impatient, then a decrease in network effects causes the online firm to decrease the advertisements and referrals in initial stages close to $t=0$ as the marginal costs for increasing the advertisements and referrals are high (at small time instances there are many sponsored arrivals which causes the marginal costs to be high). When the online firm is close to the steady state  it does not want to lose the user base acquired up until that point, and it can prevent the users from exiting  as quickly by increasing the advertisements and referrals as the marginal costs for increasing advertisements are low (at large time instances there are few sponsored arrivals which causes the marginal costs to be low).  We summarize the above comparisons in Table III. 
  
  \begin{table}
  	\centering
  	\begin{tabular}{ |l| l | l|  }
  		\hline
  		Firm type (based on $\rho$) & Change in referrals and ads &  Time \\ \hline
  		Very patient & Decrease & Steady state\\
  		\hline
  		Very impatient &  Increase & Steady state  \\ \hline
  		Very patient & Increase & Start $t=0$\\
  		\hline
  		Very impatient &  Decrease & Start $t=0$  \\ 
  		\hline
  	\end{tabular}
  	\caption{Impact of the decrease in the network effects on the referrals and advertisements in steady state}
  \end{table}

   
  Until now we have discussed the impact of the online firm's patience level on the change in the advertisements and referrals in the steady state. Next, we show that the type of the online firm as determined by the heterogeneity in its user base  also plays a key role in its decisions. In the comparisons that follow we will focus on the impact of this heterogeneity on the online firm's decisions to change the optimal policy at large values of time instances (in steady state) when the network effects decrease, while keeping other parameters of the problem fixed.

  We partition  online firms into two categories based on the rate at which the marginal benefits increase when the user base decreases. Note that $b^{''}(x)x$ is proportional to the change in the marginal benefits and $b^{'}(x)$ is proportional to the marginal benefits. Therefore, $\frac{b^{''}(x)x}{b^{'}(x)}$ is proportional the rate of change of marginal benefit. 
 We also define  constants $\Theta$, $\Delta$ such that $\Theta \leq \Delta$, whose expressions can be found in Appendix F.  
 \begin{itemize}
 	\item If  $\sup_{x \in [x_{initial},x_{u}]}\frac{b^{''}(x)x}{b^{'}(x)} < \Theta$, then the online firm is of the heterogeneous type for which the marginal benefits from adding users increase quickly when the user base decreases. These online firms have a very heterogeneous user base, for e.g.,  online  gaming firms. 
 	\item If  $\inf_{x \in [x_{initial},x_{u}]}\frac{b^{''}(x)x}{b^{'}(x)} > \Delta$, then the online firm is of the homogeneous type for which the marginal benefits from adding users increase slowly when the user base decreases. These online firms have a very homogeneous user base, for e.g., subscription based online firms.
 \end{itemize} 
%

 \begin{theorem}
 	If there is a decrease in the strength of the network effects $\eta$ and 
 	\begin{itemize}
 		\item if the online firm of the heterogeneous type, i.e.  $\sup_{x \in [x_{initial},x_{u}]}\frac{b^{''}(x)x}{b^{'}(x)} < \Theta$, then the optimal policy prescribes  to increase advertisements and referrals  in the steady state;
 		\item if the online firm of the homogenoeus type, i.e.    $\inf_{x \in [x_{initial},x_{u}]}\frac{b^{''}(x)x}{b^{'}(x)} > \Delta$, then the optimal policy prescribes  to decrease advertisements and referrals in the steady state.
 	\end{itemize}
 \end{theorem}
 
 For the proof of above theorem refer to Appendix F.

 The main conclusion from Theorem 7 is that the type of the online firm based on its  benefit function determines its response to the change in network effects. The intuition for the above theorem is the following. Suppose an online firm of the heterogeneous type  sees a decrease in network effects, which can lead to reduction in the user base. A reduction in the user base can cause a large increase in the marginal benefits from adding users,  because some of the users who generate large benefits for the online firm will have exited. Hence, by giving more advertisements and referrals the online firm can avoid a sharp reduction in the user base and thus sustain the users who generate large benefits for a longer duration. On the other hand, if the online firm is of the homogeneous type then the increase in the marginal benefits from reduction in the user base is not sufficiently high.  Hence, the online firm will save on costs by reducing advertisements and referrals.

 We summarize this comparison in Table IV.  
 \begin{table}
 	\centering
 	\begin{tabular}{ |l| l | l|  }
 		\hline
 		Firm type (speed of saturation of marginal benefits) & Change in &  Time \\ 
 		 &  referrals and ads&  \\ \hline 
 		Fast  $\sup_{x \in [x_{initial},x_{u}]}\frac{b^{''}(x)x}{b^{'}(x)} < \Theta$  (online gaming firms) & Increase & Steady state\\
 		\hline
 		Slow   $\inf_{x \in [x_{initial},x_{u}]}\frac{b^{''}(x)x}{b^{'}(x)} > \Delta$ (subscription based online firms) &  Decrease & Steady state  \\ 
 		\hline
 	\end{tabular}
 	\caption{Impact of the decrease in the network effects on the referrals and advertisements in steady state}
 \end{table}

In the next section we discuss an extension of the model where we incorporate uncertainty in user arrivals/exits in the online firm's policy design problem.
\section{Extensions}

\subsection{Impact of uncertainty in the user arrivals/exits on the optimal policy}

\subsubsection{Model:}
In this section we model the uncertainty in the user arrivals and exits and reformulate the online firm's policy design problem, which is now aimed at maximizing the expected long-term profits. A main source of uncertainty is the unpredictability in the preferences of the users for the online firm. Given the analytical difficulties associated with dealing with stochasticity in user dynamics, we will consider a specific case of the model proposed in the previous sections.

We assume that the direct arrivals as a function of the user base $X=x$ (user base size $X$ is a random variable whose dynamics are described later) are constant, i.e. $f(x)=\theta$ and the growth in the stay time due to increase in the user base user base is constant, i.e. $g(x)=1$ \footnote{$g(x)$ is constant for  online firms that do not rely significantly on the social interaction amongst users, and hence the stay time does not increase with the growth of the user base.}.   The random arrivals/exits in a time interval $dt$ is given by $\sigma  X dW_{t}$, where $\sigma$ is the standard deviation of a standard Brownian motion  and $dW_{t}$ is the differential of the Brownian motion process. Note that if $\sigma  X dW_{t}$ is negative then it indicates that more users are leaving, while if $\sigma  X dW_{t}$ is positive then it indicates that more users are joining. The online firm chooses the intensity of advertisements and referrals and as a result it controls the rate of  sponsored arrivals  $\Lambda(t,X)$ \footnote{In the deterministic case discussed we specified sponsored arrivals only as a function of time, i.e. $\lambda(t)$, as it was a sufficient succinct representation} at time $t$ and user base level $X$.
The user arrival process on the online firm follows  the  stochastic user base dynamic given as $dX = (\theta + \Lambda(t,X) - \frac{1}{(\tilde{\eta}+\eta)}X)dt + \sigma X dW_{t}$.

We analyze the benefits and costs of the online firm. We assume that the online firm has a finite capacity $\Gamma$, which is true in  practice because the servers of the online firm cannot handle unlimited traffic. The benefit function is $b:[0,\infty)\rightarrow [0,\infty)$ with $b(x) =  \Gamma^2 - (\Gamma-x)^2$. $b(x)$  increases as the user base size increases from $x=0$ up to the capacity $x=\Gamma$ beyond which it decreases and $b(x)$ is a concave function. In the previous sections we had assumed that $b(x)$ is an increasing function. Here we will assume that the capacity is sufficiently large, i.e. $\Gamma \geq \theta \eta$, which will ensure that along the optimal user base trajectory  there is typically an increase in the benefits. The long-term average benefit is given as $B(\bar{X}|\bar{\Lambda})= \int_{0}^{\infty}b(x(t))e^{-\rho t}dt$, where $\bar{X}=\{x(t)\}_{t=0}^{\infty}$ represents a sample path of the user base trajectory for a given policy $\bar{\Lambda}=\{\Lambda(t,x(t))\}_{t=0}^{\infty}$. The expectation of the long-term benefit is given as $E_{\bar{X}}[B(\bar{X}|\bar{\Lambda})]$, where the expectation is taken over the sample paths $\bar{X}$.  We assume that the flow cost  to the online firm is quadratic and is given as $c(\Lambda)=c.\Lambda^2$ \footnote{ In \cite{gourio2014customer} the authors assumed quadratic costs for searching for the customers and their calibration of the model showed that quadratic costs are a good fit for the data}. The long-term average cost is given as 
$C(\bar{\Lambda})=\int_{0}^{\infty}c(\Lambda(t,x(t)))e^{-\rho t}dt$ and the expected cost is given as $E_{\bar{X}}[C(\bar{\Lambda})]$, where the expectation is taken over the sample paths $\bar{X}$. The long-term expected profit of the online firm is given as $E_{\bar{X}}[B(\bar{X}|\bar{\Lambda})-C(\bar{\Lambda})]$. 

We state the online firm's stochastic policy design problem\footnote{In the policy design problem we do not have a constraint on the policy $\Lambda(t,X)$ being non-negative. However, if $\Gamma$ is very large then $\Lambda(t)$ will have a low probability of being less than zero. } that takes the initial user base $X(0)=x_{initial}$ as given as follows: 
\vspace{-1em}
\begin{eqnarray*}
\Pi_{stoch}(x_{initial},\eta) = && \max_{\Lambda(t,X)\in \mathbb{R}} E_{\bar{X}}[B(\bar{X}|\bar{\Lambda})-C(\bar{\Lambda})] \\
&& \text{subject to} \; dX = (\theta + \Lambda(t,X) - \frac{1}{(\tilde{\eta}+\eta)}X)dt + \sigma X dW_{t}
	\end{eqnarray*}

	We first make some remarks about our approach to solve the above problem.  $\Pi_{stoch}(x_{initial},\eta)$ has to satisfy the HJB equation for the above problem.  We formulate the HJB equation for the above problem and arrive at a solution to the HJB equation. Note that the above problem is an infinite horizon Linear Quadratic Gaussian (LQG) control problem. From \cite{speyer2008stochastic} we can conclude that the negative definite solution to the HJB equation will correspond to $\Pi_{stoch}(x_{initial},\eta)$. The optimal policy can be derived using $\Pi_{stoch}(x_{initial},\eta)$ as it is the maximizer of the term in the RHS of the HJB equation. Let us denote the  policy that we construct   as $\Lambda^{*}(t,X)$. This policy $\Lambda^{*}(t,X)$  will lead to a user base trajectory,  which is random due to the stochastic nature of the dynamic. A sample path of the optimal user base trajectory is denoted as $X_{x_{initial}}^{*}(t)$ and the corresponding sample path for the optimal policy is given as $\Lambda_{x_{initial}}^{*}(t)$, which satisfies $\Lambda_{x_{initial}}^{*}(t)=\Lambda^{*}(t,X_{x_{initial}}^{*}(t))$. In the previous sections we analyzed the deterministic counterpart of $\Lambda_{x_{initial}}^{*}(t)$ given as $\lambda^{*}_{x_{initial}}(t)$. In this section we will analyze the expectation of sponsored arrivals $\Lambda_{x_{initial}}^{*}(t)$, where the expectation is taken over the optimal sample paths  $X_{x_{initial}}^{*}(t)$ and is denoted as $E_{X_{x_{initial}}^{*}}[\Lambda_{x_{initial}}^{*}(t)]=E_{X_{x_{initial}}^{*}}[\Lambda^{*}(t,X_{x_{initial}}^{*}(t))]$. Also, define the expectation of the optimal user base trajectory  as $E_{X_{x_{initial}}^{*}}[X_{x_{initial}}^{*}(t)]$
	

We will first briefly discuss the results along the lines of the previous sections. For these results we will assume that the standard deviation of the Brownian motion is bounded above by $\rho$.   The details of the discussion that follows have been worked out in Appendix I. We can prove that the expectation of the optimal policy $E[\Lambda_{x_{initial}}^{*}(t)]$ decreases and converges to a steady state value as $t\rightarrow \infty$ (extension of Theorem 2 and 3), which shows that even though the optimal policy $\Lambda_{x_{initial}}^{*}(t)$ can take random values it will have a stationary (time-independent) mean as time goes to infinity. In Section 6 we discussed the impact of network effects on the optimal policy at the same size of the user base $\zeta^{*}(x)$ and the optimal policy at the same time $\lambda_{x_{initial}}^{*}(t)$.  If the network effects decrease, then the change in the optimal policy  depends on the capacity of the online firm $\Gamma$. See Appendix I for the details. 

Next, we discuss the impact of increase in the uncertainty in the user base $\sigma$. We first compare the optimal policy $\Lambda^{*}(t,X)$ when the network effects decrease at the same size of the user base.

\subsubsection{Impact of increase in uncertainty: comparing the optimal policy at the same size of the user base:}

  In the next theorem we show how  $\Lambda^{*}(t,X)$ changes when the level of uncertainty $\sigma$ increases. The optimal policy $\Lambda^{*}(t,X)$ does not depend on time $t$ and is in fact an affine decreasing function in $X$ (see  Appendix I for details). Hence, for a fixed $X$ the $\Lambda^{*}(t,X)$ is deterministic. 

\begin{theorem}
	If the level of the uncertainty $\sigma$ increases, then the optimal policy prescribes to decrease the referrals and advertisements for the same size of the user base $X$.
	\end{theorem}
	For the proof of the above theorem refer to Appendix I. 
	
	The intuition behind the above theorem is as follows. When  $\sigma$ increases then there is a increase in the number of random arrivals/exits. There will be an increase in benefits due to the large number of random arrivals but there will also be an increase in the losses due to the large number of random exits. Hence, the optimal policy needs to weigh these potential benefits against the potential losses and then prescribe accordingly. If there is a large number of exits, then the size of the user base decreases by a large value and it will take a large amount of time for the user base to increase in size. Hence, if the online firm increases the advertisements and referrals, then it will have to bear large costs when the user base size decreases by a large value. Since the costs are large and the  increase in the benefits due to the increase in the random user arrivals will not be sufficient to compensate for it (due to concavity of the benefit function), the online firm is prescribed to decrease advertisements and referrals. 
	
%

	
\subsubsection{Impact of increase in uncertainty: comparing the user base size and profits achieved}

In this section we aim to compare the expected size of the user base and the profits that are achieved as $t\rightarrow \infty$, i.e. in steady state, under varying levels of uncertainty. The uncertainty in the user dynamics has both a positive component in the form of random arrivals and a negative component in the form of random exits. In the next theorem we show that  under higher levels of uncertainty  in user arrivals/exits following  the optimal policy  leads to a smaller expected size of the user base in steady state as well as the expected profits achieved by the online firm. 

\begin{theorem}
	If the level of uncertainty $\sigma$ increases, then the prescriptions of the optimal policy lead to a smaller expected size of the user base and lesser expected profits that are achieved. 
	\end{theorem}
	
	For the proof of above theorem refer to Appendix I. 
		Due to an increase in the level of uncertainty  there is an increase in both random arrivals and exits. We have already seen in Theorem 7 that the optimal policy prescribes the online firm to decrease advertisements and referrals.  The increase in random arrivals cannot compensate for the increase in random exits and decrease in sponsored arrivals. Hence, the expected size of the user base that is achieved in steady state is smaller and the same reasoning applies for the expected profits as well.
\section{Conclusion}

In this work we developed a first general model to understand the complementary forces of network effects and marketing. Our key findings show how an online firm should change its optimal policy when there is a decline in network effects and we identify the main forces that determine the change in the optimal policy. If there is a decline in the network effects, then an online firm that is very patient is prescribed to increase the advertisements and referrals when the user base is small and decrease the advertisements  when the user base is large.   
On the other hand, if an online firm is very impatient, then it should decrease the advertisements and referrals when the user base is small and increase the advertisements and referrals when the user base is large. 

We also identify another interesting dimension that differentiates the response of the online firms when there is a decline in network effects. If there is a decline in the network effects and if there is significant heterogeneity in the revenue generated by the users for the online firm, then the optimal policy prescribes to increase the advertisements and referrals. On the other hand, if there is a decline in the network effects and if there is homogeneity in the revenue generated by the users for the online firm, then the optimal policy prescribes to decrease the advertisements and referrals. We also analyze the impact of uncertainty in the user base on the online firm's optimal policy,  user base size achieved and the profits achieved. If the level of uncertainty is higher, then the optimal policy prescribes the online firm to be less aggressive, and thus we can show that the expected steady state user base and the expected profits that are achieved are lower.

In this work we made the first step towards analyzing the joint forces of marketing and network effects. There are several interesting future research topics that can be analyzed. First, we can consider the setting where multiple firms compete to  build a user base. This will help understand the interaction of network effects, marketing and competition.  We considered the setting where the online firm does not influence the type distribution of the users who enter  e.g., through targeted advertisements. It will be interesting to extend the model to incorporate this dimension as well.

%
%
%

\vspace{2em}

 \begin{APPENDICES}

 \section{Proof of Existence and Uniqueness of the Optimal Policy as a function of user base and as a function of time}
 
 In Proposition 1 and 2 we state that the optimal policy  as a function of user base and the optimal policy as a function of time exist and is unique. In this Appendix we develop the proofs for both Proposition 1 and 2. Next, we outline the steps that will lead to the proof of both the propositions.

 We will first begin by showing that the optimal policy as a function of time $\lambda_{x_{initial}}^{*}$ and the corresponding user base trajectory $x_{x_{initial}}^{*}$  exist. In order to do so we will use the sufficient conditions proved in \cite{existence}. Once we show that the optimal policy as a function of time and the corresponding user base trajectory exist we will show that the value function starting from $x_{initial}$ given as $\Pi^{*}(x_{initial}, \eta)$ is finite. We can extend this and show that $\Pi^{*}(x, \eta)$ is finite for all values in a finite interval $[x_{initial},x_{u}]$.  $x_{u}= \max\{(f^{sup}+\lambda^{sup})(\tilde{\eta}+\eta^{sup} g^{sup}), x_{m} \}$, where $x_{m}$ is a constant defined in Appendix E. The reason why we define $x_{u}$ in this way will become clear later.  Next, we will show that the value function  $\Pi^{*}(x, \eta)$ is differentiable in the interval $[x_{initial},x_{u}]$. If the value function is differentiable then it will imply that it has to satisfy the HJB equation \cite{kushner1971introduction}.
 
 $\rho \Pi(x,\eta) =  \max_{0 \leq \lambda \leq \lambda^{sup}}\big(b(x)) -c(\lambda) + \frac{\partial \Pi(x,\eta)}{\partial x}(f(x)+\lambda- \frac{1}{\tilde{\eta}+ \eta g(x)}x)\big) $
 
     The maximizer of the RHS in the HJB equation corresponds to the optimal policy as a function of user base $\zeta^{*}$. We will show that the maximizer of the  HJB equation exists and is unique. This will complete the proof of Proposition 1.   Next, we will need to show that the optimal policy and the corresponding user base trajectory as a function of time is unique. For this part, we will use results from the theory of ODEs \cite{Levinson}, which will help us prove the uniqueness of the solution to the differential equations.
  
  We begin by showing that the $(\lambda_{x_{initial}}^{*}, x_{x_{initial}}^{*})$ exists.
  
   We now show that the assumptions stated in \cite{existence} are satisfied by our model.
  
  \begin{enumerate}
  	
  	\item First, the policies $\lambda(t)$ are required to be Lebesgue measurable and the user base trajectories $x(t)$ is required to be absolutely continuous on every interval $[0,T]$. Both these requirements are true due to the modeling assumptions in Section 3. Denote the rate at which the user base changes as $\frac{dx}{dt} = k_{\eta}(x,\lambda)= 
  	f(x)  +  \lambda - \frac{1}{\tilde{\eta} +  \eta g(x)}x$. We can also write $k_{\eta}(x,\lambda) =  h_{\eta}(x) +  \lambda$.
  	
  	\item  We know that $k_{\eta}(x,\lambda)$ is continuous in $(x,\lambda)$ because $h_{\eta}(x)$ is continuous and $\lambda$ is added to it, which is linear.
  	
  	\item Next, we will show that any user base trajectory that starts at $x_{initial}$ and satisfies the differential equation $\frac{dx}{dt}= f(x)  + \lambda(t) -  \frac{1}{\tilde{\eta}+\eta g(x)}x$ is always bounded from both above and below by finite constants. Notice that $\frac{dx}{dt}|_{x=x_{initial}}=k_{\eta}(0,\lambda) \geq  h_{\eta}(x_{initial})>0$, which imxlies that the user base can never fall below zero. We explain this as follows. If the user base trajectory ever falls below zero then the user base trajectory will have to cross $x_{initial}$ and will have a negative slope at $x_{initial}$, which contradicts the fact that $\frac{dx}{dt}|_{x=x_{initial}}>0$.
  	 We now show that the user base can never cross $ x^{sup}=(f^{sup}+\lambda^{sup})(\tilde{\eta}+\eta^{sup} g^{sup})$. We know that the initial user base is less than $x^{sup}$. Let us assume that the user base trajectory assumes a value greater than $x^{sup}$. Since the initial user base $x_{initial}$ is less than $x^{sup}$ and the user base trajectory is continuous, which implies that the user base will have to cross $x^{sup}$ from below. This means that for the user base trajectory to cross $x^{sup}$ approaching it from below the rate of change of user base $\frac{dx}{dt}|_{x=x^{sup}}>0$ will be positive.
  	 We can get the upper bound on $\frac{dx}{dt}|_{x=x^{sup}}$
  	 $\frac{dx}{dt}|_{x=x^{sup}} =k_{\eta}(x_{sup},\lambda) \leq  f^{sup}+ \lambda^{sup}-\frac{1}{\tilde{\eta}+\eta g^{sup}}x^{sup}\leq 0$. This contradicts the fact that $\frac{dx}{dt}|_{x=x^{sup}}>0$.  Hence, the user base trajectory always stays below $x^{sup}$.  Therefore, we can claim that the user base trajectory is always contained in $S=[x_{initial},x^{sup}]$. 
  	 So, if we add another constraint to the problem, which says that the user base is restricted to be in $S$, the optimization problem will still be equivalent. We add such a constraint to bring the problem in the form required in \cite{existence}. 
  	 
  	 \item Since $x_{initial}$ is a fixed value, which implies that the set containing the initial value is compact.
   
    \item Let us define a mapping from time $[0,\infty)$ to the set which $\lambda(t)$ always lies in. The set  $U=[0,\lambda^{sup}]$ always contains $\lambda(t)$. Hence, the mapping always takes the  value $U$ independent of time and is therefore, continuous.  
    
    \item The fifth condition  in \cite{existence} is automatically satisfied since the sets to which the user base trajectory and the optimal policy belong are uniformly bounded. 
    
    \item The function $\beta(x_{initial})=0$ (defined in \cite{existence}) for our case and hence, it is obviously continuous. 
    
    \item  Define a function $\Phi(t,x,\lambda)=(c(\lambda)-b(x))e^{-\rho t}$. The functions $b(.)$ and $c(.)$ are continuous in $x$ and $\lambda$ respectively. The function $e^{-\rho t}$ is continuous in $t$ as well. The function $-b(x)e^{-\rho t}$ and $c(\lambda) e^{-\rho t}$ is continuous as well. Hence, the sum of these functions, i.e. $\Phi(t,x,\lambda)$ is continuous as well. Since $c(\lambda)$ is convex and there is no other term in $\Phi(t,x,\lambda)$ that depends on $\lambda$, which means that the function $\Phi(t,x,\lambda)$ is convex in $\lambda$ as well.

  	\item Consider the negative part of the function  $\phi(t,x,\lambda) = (c(\lambda)-b(x))e^{-\rho t}$, which is equal to $-b(x)e^{-\rho t}$. We need to show that negative part of
  	$\phi(t,x,\lambda)$ i.e. $\lim _{T^{'}\rightarrow \infty}\int_{T^{'}}^{T^{''}}-b(x(t))e^{-\rho t} dt$  goes to zero. Observe that
  	$0 \geq \int_{T^{'}}^{T^{''}}-b(x(t))e^{-\rho t} dt \geq \int_{T^{'}}^{\infty}-b(x^{sup})e^{-\rho t} dt $
  	$=-b(x^{sup}) e^{-\rho T^{'}}\frac{1}{\rho} $. Hence, clearly the limit of the term goes to $0$. Note that this limit will be $0$ for all user base trajectories. Also, note that here we are using the property that both the Lebesgue and Reimann integral take the same value \cite{Apostol} 
  	
  	\item In addition we also know that there is a feasible trajectory i.e. which satisfies all the conditions above and gives a finite value of the objective. Consider the case when $\lambda(t)=0$ then we need to show that there is a corresponding feasible user base trajectory, which starts at $x_{initial}$ and satisfies $\frac{dx}{dt} =  h_{\eta}(x)$.  We need to show that there exists a user base trajectory over the time interval $[0,\infty)$ that satisfies the above differential equation at all times. We will show that the maximal interval of existence of the above differential equation is $[0,\infty)$.
  	 We know that $ h_{\eta}(x)$ is continuously differentiable, which also implies that it is Locally Lipschitz. Define a set $Z=(0,x^{sup}+\iota)$, where $\iota>0$. We know that the function $h_{\eta}(x)$ is bounded on $Z$ (because closure of $Z$ is compact and $h_{\eta}(x)$ is continuous on closure of $Z$). From  \cite{kartsatosadvanced} we know that maximal interval of existence is either $[0,\infty)$ or $[0,w]$ and the user base trajectory at time $w$ has to be on the boundary of $Z$. This means that the user base has to be either $0$ or $x^{sup}$. But both of these values are not possible. $x^{sup}+\iota$ and $0$ cannot be attained, the reasoning is same as in point 3 given above, where it was explained as to how every user base trajectory is bounded between $[x_{initial},x^{sup}]$. This implies that the only other possibility for maximal interval of existence is $[0,\infty)$. This shows that there is a feasible user base trajectory that satisfies the differential equation at all times. Next, we want to show that the $\int_{0}^{\infty} (b(x(t))-c(\lambda(t)))e^{-\rho t} dt$ computed using this  user base trajectory and $\lambda(t)=0$ exists. The integrand  $(b(x(t))-c(\lambda(t)))e^{-\rho t}< |b(x(t))e^{-\rho t}| < b(x^{sup})e^{-\rho t}$ and $\int_{0}^{\infty }b(x^{sup})e^{-\rho t}$ exists. Hence, the above integral exists and is finite.  
  	 	
  \end{enumerate}
  
  Given that all the assumptions required in \cite{existence} hold, we can use their main theorem, which proves the existence of an optimal policy and corresponding user base trajectory pair $(\lambda_{x_{initial}}^{*}(.), x_{x_{initial}}^{*}(.))$ exists. Let the corresponding optimal value of the objective be given as $\Pi^{*}(x_{initial},\eta)$. Note that $|\Pi^{*}(x_{initial},\eta)|$ is bounded by a constant, which can be shown using the fact that the user base trajectory is always bounded. Formally stated $|\Pi^{*}(x_{initial},\eta)|< \int_{0}^{\infty} b(x^{sup})e^{-\rho t} dt = \frac{b(x^{sup})}{\rho}$.  Similarly, we can get a lower bound $\Pi^{*}(x_{initial},\eta) \geq \frac{-c(\lambda^{sup})}{\rho}$. 
  
  Suppose we take another value for the initial user base $x^{'}_{i}$, which is between $[x_{initial},x_{u}]$ then also the above proof of existence will continue to hold. This means that the objective's optimal value is well defined and finite for all initial user base values between $[x_{initial},x_{u}]$. We will only need to make some changes in the approach to satisfy condition 3, where we show that the user base trajectory is bounded. We showed in condition 3 that the user base is bounded above by $x^{sup}$. Now since the initial user base can be more than $x^{sup}$ the bound will change to $x_{u}$ because if the initial user base is more than $x^{sup}$ then the user base trajectory will always decrease. Hence, in this case the user base trajectory is in the interval $[x_{initial},x_{u}]$. In this way we can extend the domain  of $\Pi^{*}(x,\eta)$ to $[x_{initial},x_{u}]$. In fact observe that the domain can be extended to any interval $[0,s]$, where $s$ is a finite constant, because the derivation did not rely on the fixed values $x_{initial}$, $x_{u}$.

   Note that we will in fact only take a fixed value of initial user base for the entire paper $x_{initial}$. The above definition of the value function on the entire interval is to consider the optimal policy at different user base levels that will be achieved by the user base trajectory from initial user base upto the steady state. 
  
  Next, we want to show that the value function $\Pi^{*}(x,\eta)$ is differentiable in the interior of $[x_{initial},x_{u}]$. We can define another set $S_{in}=[x_{initial}-\iota,x_{u}+\iota]$, where $0<\iota<x_{initial}$. The user base trajectory starting at $x_{initial}$ is always in the interior of $S_{in}$. We know that $\frac{dx}{dt}= h_{\eta}(x) + \lambda \leq f^{sup}+\lambda^{sup}$. The lower bound is given as $h_{\eta}(x) + \lambda \geq f^{inf}- \frac{1}{\tilde{\eta}}x_{u}$. These upper and lower bounds are derived based on the following facts about the model, the direct arrival rate is bounded above by and below by $f^{sup}$ and $f^{inf}$, the sponsored arrival rate is bounded above and below $\lambda^{sup}$ and $0$, the user base is in set $S$ and the stay time of the users is bounded above and below as well. While the bounds that we arrive at are weak and can be made tighter. As we will see that we do not need tight bounds for our purposes. 
  
  Therefore, we can say that $\frac{dx}{dt}$ always lies in the interior of $\Pi_{in}=[f^{inf}- \frac{1}{\tilde{\eta}}x_{u}-\omega, f^{sup}+\lambda^{sup} +\omega]$. The user base trajectory and the corresponding rate of change $(x(t), \frac{dx}{dt})$ lies in the interior of the set $T=S_{in} \times \Pi_{in}$, where $\times$ is the Cartesian products of the sets.
 We now show that the assumptions in \cite{Benveniste} are satisfied and hence, we can prove that the optimal value function is differentiable. 
 
 \begin{enumerate}
 	\item The set $T$ is a finite interval and hence, it is convex and has a non-empty interior. 
 	
 	\item Substitute $\lambda(t) = \frac{dx}{dt} - h_{\eta}(x)$ in $\phi(t,x,\lambda)=(c(\lambda)-b(x))e^{-\rho t}$ to get $\phi(t,x,\frac{dx}{dt})=(c(\frac{dx}{dt}-h_{\eta}(x))-b(x))e^{-\rho t}$. We know that the $h_{\eta}(x)$ is concave.  Hence, we can see that $\frac{dx}{dt}-h_{\eta}(x)$ is jointly convex in $\frac{dx}{dt}$ and $x$. We know that $c(\lambda)$ is increasing and convex. But note that the domain of $c(.)$ as assumed in our model is only the non-negative real axis.
 	We want to show that $c(\frac{dx}{dt}-h_{\eta}(x))$
  But in order to be able to use composition rules, which use sufficient conditions to show the convexity of   $c(\frac{dx}{dt}-h_{\eta}(x))$ we require that the $c$ is defined over the entire real line. We can extend the definition of $c(.)$ and assume that it is also defined on the negative half of the real line. Even under the extended definition $c(.)$ is strictly convex and increasing. Note that having this extension does not affect the policy design problem because we have a hard constraint that $\lambda$ is non-negative. So, now we can state that $c(\frac{dx}{dt}-h_{\eta}(x))$ is convex. Hence,  the function $\phi(t,x,\frac{dx}{dt})$ is convex in $(x, \frac{dx}{dt})$. 
  
  \item  The optimal value function $\Pi^{*}(x,\eta)$ is well defined and bounded for all values in the interval $S_{in}$. For every value in the interior of the interval $S_{in}$ the value function is well defined in its neighborhood. 
  
  \item The optimal user base trajectory and the corresponding $\frac{dx}{dt}$ always lie in the interior of $T$. Since the trajectories are always in the interior we can deduce that there is an $\epsilon$ ball around these trajectories, which is contained in $T$. 
  
  Note that the above four conditions are sufficient to show that $\Pi^{*}(x,\eta)$ is differentiable. We also get that $\Pi^{*}(x,\eta)$ is continuously differentiable. Now we want to obtain the optimal policy as a function of user base $\zeta^{*}(.)$. We know that the $\zeta^{*}(x)$ corresponds to $ \arg \max_{0 \leq \lambda \leq \lambda^{sup}}\big(b(x)) -c(\lambda) + \frac{\partial \Pi(x,\eta)}{\partial x}(f(x)+\lambda- \frac{1}{\tilde{\eta}+ \eta g(x)}x)\big)$ (due to Bellman's principle of optimality). This maximization is a simple one variable constrained optimization problem. Note that there are two constraints $\lambda\leq \lambda^{sup}$ and $-\lambda\leq 0$, we introduce Lagrange multipliers corresponding to these constraints $m_{1}$ and $m_{2}$. We now state the  KKT conditions for the above problem. The above problem is concave and the Slater's conditions (there is a non empty interior for the constraint set) are satisfied, hence, KKT conditions given below are necessary and sufficient (\citealt{Boyd}). 
  
 \begin{eqnarray*}
 &&	m_{1}\geq 0,\; m_{2}\geq 0 ]\; \\ 
 &&	\lambda \leq \lambda^{sup},\; -\lambda \leq 0\; \\
 &&	m_{1}(\lambda-\lambda^{sup}) =0,\; m_{2}(\lambda)=0\;  \\
  &&  c^{'}(\lambda)- \frac{\partial \Pi^{*}(x,\eta)}{\partial x} + m_{1} -m_{2}=0
\end{eqnarray*}
  
  We look at the solution to the above problem in cases. 
  \begin{enumerate}
  	\item Consider the case when $\frac{\partial \Pi^{*}(x,\eta)}{\partial x}>c^{'}(\lambda^{sup})$. In this case the first possibility is that $\lambda=\lambda^{sup}$ is a solution. If that is the case then $m_{2}=0$  and $m_{1}=\frac{\partial \Pi^{*}(x,\eta)}{\partial x}-c^{'}(\lambda^{sup})>0$. Hence, in this case $\lambda^{sup}$, $m_{2}=0$ and $m_{1}=\frac{\partial \Pi^{*}(x,\eta)}{\partial x}-c^{'}(\lambda^{sup})$ correspond to the unique solution of the optimization problem.
  	
  	\item Consider the case when $\frac{\partial \Pi^{*}(x,\eta)}{\partial x}<c^{'}(0)$. In this case the first possibility is that $\lambda=\lambda^{sup}$ is a solution. If that is the case then $m_{1}=0$  and $m_{2}=c^{'}(0)-\frac{\partial \Pi^{*}(x,\eta)}{\partial x}>0$.  Hence, in this case $\lambda=0$, $m_{1}=0$ and $m_{2}=c^{'}(0)-\frac{\partial \Pi^{*}(x,\eta)}{\partial x}$ correspond to the unique solution of the optimization problem.
  	
  	\item Consider the case when $ c^{'}(0)\leq \frac{\partial \Pi^{*}(x,\eta)}{\partial x}\leq c^{'}(\lambda^{sup})$.  From complementary slackness we can see that $m_{1}$ and $m_{2}$ will be zero. In this case $c^{'}(\lambda)=\frac{\partial \Pi(x,\eta)}{\partial x}$ and in this case we can use the property that $c^{'}(\lambda)$ is strictly increasing (due to strict convexity of $c(.)$) to conclude that the solution to $\lambda$ in this case will be unique and we can concisely write it as $\lambda= c'^{-1}(\frac{\partial \Pi(x,\eta)}{\partial x})$. Since $ c^{'}(0)\leq \frac{\partial \Pi^{*}(x,\eta)}{\partial p}\leq c^{'}(\lambda^{sup})$ we can see that $\lambda= c'^{-1}(\frac{\partial \Pi(x,\eta)}{\partial x})$ is primal feasible. Hence, in this case $\lambda=c'^{-1}(\frac{d\Pi(x,\eta)}{dx})$, $m_{1}=0$ and $m_{2}=0$ correspond to the unique solution of the optimization problem.
  	
 From the above discussion we can conclude that the optimizer as a function of the user base exists and is unique. In proving the uniqueness above we did not use the Assumption 1 part 1. In fact we require Assumption 1 part 1 for the uniqueness of the optimal policy as a function of time. 
 
  	So, up until now we have shown that the optimal policy as a function of user base exists and is unique, which proves Proposition 1. We have also proven part of Proposition 2, since we know that the optimal policy and the corresponding user base trajectory exists. Next, we aim to show that the optimal policy as a function of time is unique. The optimal user base trajectory satisfies the following differential equation, $\frac{dx}{dt}=  f(x)+\lambda(t)- \frac{1}{\tilde{\eta}+\eta g(x)}x$. We already know that the optimal policy as a function of user base is unique and is given as $\zeta^{*}(.)$. If the user base at time $t$ is $x(t)$ then the optimal policy at time $t$ is  $\zeta^{*}(x(t))$ (This is true because of the following two reasons, i) optimal policy as a function of time $\lambda^{*}_{x_{initial}}$ needs to satisfy Bellman's principle of optimality, ii) optimization problem has an infinite horizon). Hence, we can formulate the differential equation that the optimal user base trajectory needs to satisfy is given as  $\frac{dx}{dt}=  f(x)+\zeta^{*}(x)- \frac{1}{\tilde{\eta}+\eta g(x)}x$. If we can show that the solution to this differential equation is unique then we are set because then $\zeta^{*}(x(t))$ will be uniquely determined. We can write the differential equation as $\frac{dx}{dt}=  h_{\eta}(x)+\zeta^{*}(x)$. The first term $h_{\eta}(x)$ is continuously differentiable on $(x_{initial},x_{u})$ and if the second term is continuously differentiable as well then $h_{\eta}(x)+\zeta^{*}(x)$ will be continuously differentiable on $(x_{initial},x_{u})$.  We know that an optimal user base trajectory exists defined over the time interval $[0,\infty)$. Thus, we can use standard results on uniqueness in ODEs, which extend Picard's results to global uniqueness \cite{Levinson}. Next, we still need to ensure that the second term is differentiable.

   $ c^{'}(0)\leq \frac{\partial \Pi^{*}(x,\eta)}{\partial x}\leq c^{'}(\lambda^{sup})$ holds for all user base levels in the interval $[x_{initial},x_{u}]$. We know that the user base trajectory starting at $x_{initial}<f^{inf}\tilde{\eta}$ always lies in the interval $[x_{initial},x_{u}]$ and thus, we can conclude the sponsored arrivals at all times will be given by 
  $c'^{-1}(\frac{\partial \Pi^{*}(x,\eta)}{\partial x})$. This condition ensures that the optimal policy is always in the interior of $ [0,\lambda^{sup}]$.  If $\frac{\partial \Pi^{*}(x,\eta)}{\partial x}$ is continuously differentiable then we can conclude that $c'^{-1}(\frac{\partial \Pi^{*}(x,\eta)}{\partial x})$ is continuously differentiable. 
  
 
  \end{enumerate}
  
  Besides the above cases there can be scenarios when the above conditions do not hold. In those cases it may be difficult to show the uniqueness of the optimal policy as a function of time. In the proofs of theorems that will follow we will discuss the impact of the case when there is non-uniqueness. Essentially, we will argue that all the results that we will present (Theorems 4-6) do not rely on uniqueness of the optimal policy. The results presented in those theorems will hold for all the optimal policies, if there are multiple of them.

 	\end{enumerate}  
 	
\section{Proof of Existence and Uniqueness of  the Steady state}  

In Theorem 1 we state that for every online firm in our model there is a steady state, which is unique. We first show the existence part of the proof.  In steady state we know that the user base is fixed, which means $\frac{dx}{dt}=0$. Hence, $\lambda=-h_{\eta}(x)$ is the first condition for the steady state. In Proposition 2 we used the following fact that $ c^{'}(0)< \frac{\partial \Pi^{*}(x,\eta)}{\partial x} < c^{'}(\lambda^{sup}) $ for all $ x \in [x_{initial},x_{u}]$. If this is true then the optimal policy is always given as $c'^{-1}(\frac{\partial \Pi^{*}(x,\eta)}{\partial x})$. Hence, we can write the rate of change of user base as 
$\frac{dx}{dt}=h_{\eta}(x) + c'^{-1}(\frac{\partial \Pi^{*}(x,\eta)}{\partial x})$. 
  $h_{\eta}(x_{initial}) + c'^{-1}(\frac{\partial \Pi^{*}(x_{initial},\eta)}{\partial x})$ is greater than zero. Also, note that $h_{\eta}(x_{u}) + c'^{-1}(\frac{\partial \Pi^{*}(x_{u},\eta)}{\partial x})\leq h_{\eta}(x_{u}) + \lambda^{sup} \leq 0$. Therefore, either $h_{\eta}(x_{u}) + c'^{-1}(\frac{\partial \Pi^{*}(x_{u},\eta)}{\partial x})\leq h_{\eta}(x_{u}) + \lambda^{sup} = 0$, which proves that $x_{u}$ is steady state or  $h_{\eta}(x_{u}) + c'^{-1}(\frac{\partial \Pi^{*}(x_{u},\eta)}{\partial x})\leq h_{\eta}(x_{u}) + \lambda^{sup} <0$. If the latter is true then from continuity of $h_{\eta}(x_{initial}) + c'^{-1}(\frac{\partial \Pi^{*}(x_{initial},\eta)}{\partial x})$ we know that there exists a steady state. 

 Suppose that there are two steady states $x_{s}^{1}$ and $x_{s}^{2}$. The corresponding values of the policies are as follows $\lambda_{s}^{1}=-h_{\eta}(x_{s}^{1})>0$ and $\lambda_{s}^{2}=-h_{\eta}(x_{s}^{2})>0$. Without loss of generality consider $x_{s}^{1}<  x_{s}^{2}$. $\lambda_{s}^{1}> 0$ which gives $h_{\eta}(x_{s}^{1})<0$ and similarly $h(x_{s}^{2})<0$. Recall that $h_{\eta}(x_{initial})>0$  and $h_{\eta}(x)$ is concave, which we can use to deduce that $h_{\eta}(x_{s}^{1})>h_{\eta}(x_{s}^{2})$, which implies $\lambda_{s}^{1}<\lambda_{s}^{2}$. 
 We know that $\frac{\partial \Pi^{*}(x_{s}^{1},\eta)}{\partial x} = c^{'}(\lambda_{s}^{1})$  and 
 $\frac{\partial \Pi^{*}(x_{s}^{2},\eta)}{\partial x} = c^{'}(\lambda_{s}^{2})$. Note that $c^{'}$ is an increasing function and $\frac{\partial \Pi^{*}(x_{s}^{1},\eta)}{\partial x}$ is a decreasing function in $x$ (due to concavity Lemma 1). RHS of steady state equation corresponding to $x_{s}^{2}$ ($\frac{\partial \Pi^{*}(x_{s}^{2},\eta)}{\partial x} = c^{'}(\lambda_{s}^{2})$) is higher than the RHS of steady state of $x_{s}^{1}$, while its LHS is lesser. Thus, the second equality cannot hold, which leads to a contradiction. 

\section{Proof of Convergence to the Steady state}

In the previous section we proved that there exists a steady state. In this section we prove Theorem 2. We know that there exists a unique user base trajectory that satisfies the differential equation $\frac{dx}{dt} = h_{\eta}(x) + c'^{-1}(\frac{\partial \Pi^{*}(x,\eta)}{\partial x})$ starting at $x_{initial}$. Also, we know that $h_{\eta}(x) + c'^{-1}(\frac{\partial \Pi^{*}(x,\eta)}{\partial x})$ is positive in the interior of the interval $[x_{initial},x_{s}]$. And if we choose a user base level sufficiently close $\epsilon$ to $x_{s}$ and it is given as $x_{s}-\epsilon$. Then $h_{\eta}(x) + c'^{-1}(\frac{\partial \Pi^{*}(x,\eta)}{\partial x})\geq \epsilon^{'}$ for all the values of $x \in [x_{initial}, x_{s}-\epsilon]$. This means that $\frac{dx}{dt}\geq \epsilon^{'}$ up until $x > x_{s}-\epsilon$. This means that starting with $x_{initial}$ the user base trajectory will reach the threshold in at most $\frac{x_{s}-\epsilon}{\epsilon^{'}}$ time, which is finite. In this way by repeating such an argument we can show that the user base trajectory is bound to converge.  As soon as the $x(t)$ reaches $x_{s}$ the rate of change will be zero. The corresponding optimal policy will also converge as the user base trajectory converges.

\section{Proof of Monotonicity of the Optimal policy }

	We state a lemma below, which states that the optimal value function is concave in $x$.  We discuss this Lemma separately as we need it even at later stages in other proofs. However, the proof of the differentiability in \cite{Benveniste} also proves this lemma as well.
	\begin{lemma}
		The optimal value function $\Pi^{*}(x,\eta)$ defined on the interval $[x_{initial},x_{u}]$ is concave in $x$. 
	\end{lemma}
	We know that $b(x)-c(\frac{dx}{dt}-h_{\eta}(x))$	is jointly concave in $ (x,\frac{dx}{dt})$ and $(x,\frac{dx}{dt}) \in T$, where $T$ is convex. This establishes the concavity of $\Pi^{*}(x,\eta)$.

We already know from the previous subsection that the optimal user base trajectory converges.  Clearly since $ h_{\eta}(x) + c'^{-1}(\frac{\partial \Pi^{*}(x,\eta)}{\partial x})\geq 0$ at all points on the path before the policy converges, which implies that the user base increases strictly before converging.  The corresponding optimal policy as a function of time is $\lambda_{x_{initial}}^{*}(t)= c'^{-1}(\frac{d\Pi^{*}(x_{x_{initial}}^{*}(t),\eta)}{dx})$. Note that $c'^{-1}(.)$ is increasing and  $\frac{\partial \Pi^{*}(x,\eta)}{\partial x}$ is a decreasing function in $x$ because $\Pi^{*}(x,\eta)$ is concave. These factors combined with the fact that $x_{x_{initial}}^{*}(t)$ is an increasing function of time lead us to conclude that the optimal policy will decrease with time. 

\section{Proof of results on impact of network effects on optimal policy at the same size of the user base}
In this section where we compare the optimal policy as a function of user base $\zeta^{*}$. Up until now we have shown that the optimal policy as a function of user base is well-defined for $[x_{initial}, x_{u}]$. 

In Theorem 4 we show if the online firm is patient, i.e. a low $\rho$ then the change in the optimal policy depends on its user base size. We are interested in analyzing the case when the solution $\lambda$ lies in the interior of the interval $[0,\lambda^{sup}]$. We are only interested to analyze the impact of network effects on this case because otherwise when the solution lies on the boundary then the change in the network effects will not change the solution. When the solution lies in the interior then the HJB equation can be written as follows. Let us denote $\frac{\partial \Pi^{*}(x,\eta)}{\partial x} = \Pi'^{*}(x,\eta)$ and $c'^{-1}(\lambda)=q(\lambda)$

\begin{eqnarray*}
&& \rho \Pi^{*}(x,\eta) = b(x)-c(q(\Pi'^{*}(x,\eta))) + \Pi'^{*}(x,\eta) (h_{\eta}(x)+q(\Pi'^{*}(x,\eta))) 
\end{eqnarray*}

\begin{eqnarray*}
&&\rho\frac{\partial \Pi^{*}(x,\eta)}{\partial \eta} = -c'(q'(\Pi'^{*}(x,\eta))) \frac{\partial \Pi'^{*}(x,\eta)}{\partial \eta} + (h_{\eta}(x)+q(\Pi'^{*}(x,\eta)))\frac{d\Pi'^{*}(x,\eta)}{d\eta}\\ &&\Pi'^{*}(x,\eta)q'(\Pi'^{*}(x,\eta))\frac{\partial \Pi'^{*}(x,\eta)}{\partial \eta}  +
 \Pi'^{*}(x,\eta) \frac{\partial h_{\eta}(x)}{\partial x} 
\end{eqnarray*}
We will rearrange the above equation to bring all the terms that do not contain $\frac{d\Pi'^{*}(x,\eta)}{d\eta}$ on the LHS to obtain $\rho \frac{\partial \Pi^{*}(x,\eta)}{\partial \eta} - \Pi'^{*}(x,\eta) \frac{dh_{\eta}(x)}{dx}$.  Simplifying the LHS term
\begin{equation} 
\rho \frac{\partial \Pi^{*}(x,\eta)}{\partial \eta} - \Pi'^{*}(x,\eta) \frac{1}{(\tilde{\eta} + \eta g(x))^2}g(x)x.
\end{equation}

 Next, we will show that $\frac{\partial \Pi^{*}(x,\eta)}{\partial \eta}$ is positive. The intuition for this result is that when there is an increase in the network effects then the users stay longer at the same size of the user base, which causes the user base that is active at any time $t$ to be higher.

In order to figure out $\frac{\partial \Pi^{*}(x,\eta)}{\partial \eta}$ we will find a lower bound for it. Consider $\eta$ and $\eta +\Delta \eta$. The optimal policy as a function of user base at $\eta$ is given as $\zeta_{\eta}^{*}(.)$. Construct another policy which is given as $\zeta_{\eta+\Delta \eta} = h_{\eta}(x) -h_{\eta+\Delta \eta} + \zeta_{\eta}^{*} $. Note that $\Delta \eta$ can be chosen to be sufficiently small such that $h_{\eta} -h_{\eta+\Delta \eta} + \zeta_{\eta}^{*} $ is in the interior of the set $[0,\lambda^{sup}]$. The motivation to construct this policy is as follows. Consider the rate of change of the user base when network effects is $\eta$ and the policy is $\zeta_{\eta}^{*}$, $\frac{dx_{\eta}}{dt} = \zeta_{\eta}^{*} + h_{\eta} $. Consider the rate of change of the user base when network effects is $\eta+ \Delta \eta$ and the policy is $\zeta_{\eta+\Delta \eta}$,  $\frac{dx_{\eta +\Delta \eta}}{dt} = \zeta_{\eta+\Delta \eta} + h_{\eta+\Delta \eta}$. We simplify the rate of change $\frac{dx_{\eta +\Delta \eta}}{dt}$ by substituting for $\zeta_{\eta+\Delta \eta}$ to get $\frac{dx_{\eta +\Delta \eta}}{dt} = \zeta_{\eta}^{*} + h_{\eta}$. Now we see that the rate of change  $\frac{dx_{\eta}}{dt} = \frac{dx_{\eta+\Delta \eta}}{dt}$ and since both $x_{\eta}(t)$ and $x_{\eta+\Delta \eta}$ start at the same point. Our other assumptions ensure that $\zeta_{\eta}^{*} + h_{\eta}$ satisfies conditions for the uniqueness of the solution to the differential equation.  Hence, both $x_{\eta}(t) = x_{\eta+\Delta \eta}(t)$ and we denote it for simplicity as $x(t)$. Clearly, the benefits achieved along both the policies is the same. But the latter policy needs to give out fewer advertisements as $h_{\eta+\Delta\eta}(x)> h_{\eta}(x)$. Therefore, the cost is lower at all times, which leads to a positive change in the value function. Let us compute the difference in the costs at all times. We will use the notation $\lambda_{\eta}(t)$ and $\lambda_{\eta+\Delta \eta}(t)$ for the policies at time $t$ for simplicity.  $c(\lambda_{\eta}(t))-c(\lambda_{\eta+\Delta\eta}(t))=c(\lambda_{\eta}(t)) - c(\lambda_{\eta}(t)+h_{\eta}(x(t))-h_{\eta+\Delta \eta}(x(t)))$. Since $\Delta \eta$ is small this difference can be simplified to $c^{'}(\lambda_{\eta}(t))(h_{\eta+\Delta \eta}(x(t))-h_{\eta}(x(t)))$. Hence, we can write $\frac{dc(\lambda_{\eta}(t))}{d\eta} = c^{'}(\lambda_{\eta}(t)) \frac{dh_{\eta}(x(t))}{d\eta}= c^{'}(\lambda_{\eta}(t)) \frac{g(x(t))x(t)}{(\tilde{\eta}+\eta g(x(t)))^2}$. We know that the second policy that we obtained is not necessarily optimal, but we can use the difference in the value function of the optimal policy at $\eta$ and the suboptimal policy at $\eta+\Delta \eta$. We know that there is no difference in the long-term benefits. Let us compute the difference between the long-term costs. $C(\lambda_{\eta}(.))-C(\lambda_{\eta+\Delta \eta}(.)) =  \int_{0}^{\infty}(c(\lambda_{\eta}(t)-c(\lambda_{\eta + \Delta \eta}(t))e^{-\rho t} dt $. Using the  $\frac{dc(\lambda_{\eta}(t))}{d\eta} = c^{'}(\lambda_{\eta}(t)) \frac{dh_{\eta}(x(t))}{d\eta}= c^{'}(\lambda_{\eta}(t)) \frac{g(x(t))x(t)}{(\tilde{\eta}+\eta g(x(t)))^2}$ we can get $\frac{dC(\lambda_{\eta}(.))}{d\eta} = \int_{0}^{\infty}\frac{dc(\lambda_{\eta}(t))}{d\eta} e^{-\rho t} dt = \int_{0}^{\infty} c^{'}(\lambda_{\eta}(t)) \frac{g(x(t))x(t)}{(\tilde{\eta}+\eta g(x(t)))^2} e^{-\rho t} dt$. Since each of the terms in the integrand is strictly positive we can deduce that $\frac{dC(\lambda_{\eta}(.))}{d\eta}>0$. Hence, we get $\frac{\partial \Pi^{*}(x,\eta)}{\partial \eta} \geq \frac{dC(\lambda_{\eta}(.))}{d\eta}>0 $. Next, we will show that if $\rho$ is sufficiently small then $\rho \frac{\partial \Pi^{*}(x,\eta)}{\partial \eta} - \xi'^{*}(x,\eta) \frac{1}{(\tilde{\eta} + \eta g(x))^2}g(x)x$ is negative. 

First note that $\Pi'^{*}(x,\eta)$ is positive. The intuition for this is as follows, consider two user base levels $x$ and $x+\Delta x$, where the latter has a higher value. If the online firm starts at a higher user base level $x+\Delta x$ then even if the online firm follows the same optimal policy as it did starting at $x$ then also the user base trajectory starting at $x+\Delta x$ will be higher. A higher value for user base at all times will lead to higher  benefits and the costs will be the same. This leads to an improvement in the $\Pi^{*}(x,\eta)$. In \cite{Benveniste} the authors in fact gave a technique to simplify the
$\Pi'^{*}(x,\eta)$. In our case if we use their technique we get $\Pi'^{*}(x,\eta) = c^{'}(\zeta^{*}(x))$. Although, the simplification does not give an explicit value, it still depends on the $\zeta^{*}(x)$, which should not matter for us because as we want to arrive at sufficient conditions to make $\rho \frac{\partial \Pi^{*}(x,\eta)}{\partial \eta} - \Pi'^{*}(x,\eta) \frac{1}{(\tilde{\eta} + \eta g(x))^2}g(x)x$ negative. We can get a lower bound on $c^{'}(\zeta^{*}(x))$, which we know is $c^{inf}$.  Define $\zeta = 
\sup_{x_{initial} \leq x \leq x_{u}, 0 \leq \eta \leq \eta^{sup}, 0 \leq \rho \leq \delta} \frac{\partial \Pi^{*}(p,\eta)}{\partial \eta}$. Note that the supremum  involves three terms the user base level, network effects and the discount rate
. Also, if  $\frac{\partial \Pi^{*}(x,\eta)}{\partial \eta}$ is continuous $\eta, x, \rho$ then the supremum is finite and positive. This is because we know that $\frac{\partial \Pi^{*}(x,\eta)}{\partial \eta}$ is positive and a continuous function on a compact set is bounded \cite{Rudin}. Define the infimum of $\zeta^{'}=\inf_{ x_{initial}\leq p\leq x_{u}}\frac{1}{(\tilde{\eta} + \eta g(x))^2}g(x)x$. Now we are ready to arrive at  the bound for $\rho$. $\rho^{lb}= c^{inf}\frac{\zeta^{'}}{\zeta}$. If $\rho \leq \rho^{lb}$ and if $\rho \leq \delta$ then the term $ \rho \frac{\partial \Pi^{*}(x,\eta)}{\partial \eta} - \Pi'^{*}(x,\eta) \frac{1}{(\tilde{\eta} + \eta g(x))^2}g(x)x$ is negative. We needed the additional condition in the form of $\rho \leq \delta$ because we considered the only the $\rho$'s that satisfy this constraint while deriving $\zeta$.  Hence, we can define the constraint on $\rho$ as $\rho \leq \bar{\rho} = \min\{\rho^{lb},\delta\}$. In fact we can also define a stricter threshold $\rho\leq \rho^{l}$, where $\rho^{l}=\min\{\rho^{l}, \rho_{1}^{l}\}$. The expression for $\rho_{1}^{l}$ is given in Appendix F. Note that $\rho^{l}$ is defined in such a way to have a common lower bound defined for $\rho$.

 	The term in the RHS is  
 	\begin{equation}
 	\frac{\partial \Pi^{*}(x,\eta)}{\partial \eta}(h_{\eta}(x)+q(\Pi'^{*}(x,\eta)) + \Pi'^{*}(x,\eta)q'(\Pi'^{*}(x,\eta) -c'(q'(\Pi'^{*}(x,\eta))))
 	\end{equation}
 	 If $h_{\eta}(x_{initial}) >c'(q'(\lambda^{sup}))$
 	 then we have the condition that $(h_{\eta}(x)+q(\Pi'^{*}(x,\eta)) + \Pi'^{*}(x,\eta)q'(\Pi'^{*}(x,\eta)) -c'(q'(\Pi'^{*}(x,\eta)))) $ is greater than zero in a neighborhood close to the initial user base. Hence, $\frac{\partial \Pi^{*}(x,\eta)}{\partial \eta}<0$ for all the user base levels close to $x_{initial}$. In fact we can define a threshold $x_{1}$, where $h_{\eta}(x_{1}) =c'(q'(\lambda^{sup}))$. If $x\leq x_{1}$, then $h_{\eta}(x)<c'(q'(\lambda^{sup}))$. Therefore, we can state the same for $c'^{-1}(\frac{\partial \Pi^{*}(x,\eta)}{\partial \eta})$, which implies that the optimal policy increases when there is a decrease in the network effects at low user base levels.
 	
 	 We know that $\Pi'^{*}(x,\eta)$ is bounded above and below, because we are considering the case when the optimal policy is in the interior. Thus, $q(\Pi'^{*}(x,\eta)) + \Pi'^{*}(x,\eta)q'(\Pi'^{*}(x,\eta) -c'(q'(\Pi'^{*}(x,\eta)))$ is bounded above and below. Also, note that $h_{\eta}(x)$ is a continuous function, which approaches $-\infty$ as $x$ approaches $\infty$. This means at large enough values of $x$ $h_{\eta}(x)+q(\Pi'^{*}(x,\eta)) + \Pi'^{*}(x,\eta)q'(\Pi'^{*}(x,\eta)) -c'(q'(\Pi'^{*}(x,\eta)))$ is negative. Let us now derive the threshold on the value  of $x$ at which the expression above is negative. We know that $q(\Pi'^{*}(x,\eta))$ corresponds to the optimal policy, which is bounded above by $\lambda^{sup}$. We know that $\Pi'^{*}(x,\eta) \leq c^{'}(\lambda^{sup})$. We can arrive at the supremum of $x.q'(x)-c'(q'(x))$ in the interval $ c^{inf} \leq x \leq c^{'}(\lambda^{sup})$ and let us denote it by $q^{sup}$. We know that $f(x)$ is bounded above by $f^{sup}$ and the lower bound of $\frac{1}{\tilde{\eta}+\eta g(x)}x$ is $\frac{1}{\tilde{\eta} + \eta g^{sup}}x$. In this expression $h_{\eta}(x)+q(\Pi'^{*}(x,\eta)) + \Pi'^{*}(x,\eta)q'(\Pi'^{*}(x,\eta))-c'(q'(\Pi'^{*}(x,\eta)))$ we will replace each positive term by its upper bound and the negative term by the lower bound and thus obtain the upper bound for the expression. Then, we will derive a sufficiently large value of $x$ when the upper bound is negative.  Based on the above substitutions we get $x \geq x_{2}= (f^{sup}+q^{sup}+\lambda^{sup}) (\tilde{\eta} + \eta^{sup} g^{sup})$ then the expression is negative. We can also define $x_{m}$, $x_{m}=(f^{sup}+q^{sup}+\lambda^{sup}) (\tilde{\eta} + \eta^{sup} g^{sup})+\tilde{k}$, where $\tilde{k}>0$.
 	
 	 Hence, $\frac{d\Pi'^{*}(x,\eta)}{d\eta}>0$. This means that the online firm decreases the advertisements and referrals when there is a decline in the network effects. This proves part 1 and 2 of Theorem 4. 
 	
 	Before we prove the remaining parts of Theorem 4 we want to make a remark about the previous proof. We showed sharp results about what happens when the user base is close to initial user base level and then at very high levels. First of all the bounds that have been arrived are not tight and can be potentially be made tighter. Secondly, the part of the result which discusses the comparison at high user base levels but does not discuss whether the user base level is higher than the steady state level or lesser. If it is higher then it will not be achieved along the optimal population trajectory. But if it is lesser, for e.g., consider the scenario when the following is true $x.q'(x)\leq c'(q'(x))\; \forall x \in [x_{initial},x_{u}]$ then at user base levels close to the steady state the switch (from increasing to decreasing advertisements) in the behavior of the online firm will be observed.

 	If $\rho$ is sufficiently high then just the opposite effects of the previous theorem are observed. We already know that the RHS in (2) changes sign when user base increases from low to high values. All we need to show that if $\rho$ is sufficiently high then the sign of the LHS in (1) is always positive. We can arrive at the supremum of $ \tilde{\zeta}^{'} = \sup_{x \in [x_{initial},x_{u}]} \frac{g(x)x}{\tilde{\eta} + \eta g(x)^2}$. We arrive at the infimum of and denote it as $\tilde{\zeta} = \inf_{x\in [x_{initial},x_{u}],  0 \leq \eta \leq \eta ^{sup}, \rho \geq 0}\frac{\partial \Pi^{*}(x,\eta)}{\partial \eta}$. We know that $\tilde{\zeta}$ is non-negative, we assume that it is strictly greater than zero.  Hence, if $\rho \geq c^{'}(\lambda^{sup})\frac{\tilde{\zeta}^{'}}{\tilde{\zeta}}$ then the term $\rho \frac{\partial \Pi^{*}(x,\eta)}{\partial \eta} - \Pi'^{*}(x,\eta) \frac{dh_{\eta}(x)}{dx}$ is positive. We can define the constraint on $\rho$ as $\rho \geq \rho^{u}= \max\{\rho_{1}^{u},c^{'}(\lambda^{sup})\frac{\tilde{\zeta}^{'}}{\tilde{\zeta}}\}$ This proves part 3 and 4 of Theorem 4. 
 	

\section{Proof of comparison in the steady state}

Consider the case when $\lambda_{s}$ occurs in the interior of the interval $[0,\lambda^{sup}]$ then we know that the following condition has to hold $c^{'}(\lambda_{s})= \frac{d\Pi^{*}(x_{s},\eta)}{dx}$. Our first objective is to simplify the expression in $\frac{d\Pi^{*}(x_{s},\eta)}{dx}$. Consider a function $W(x,\eta)= \frac{b(x)-c(-h(x))}{\rho}$. $W(x,\eta)$ is concave for the reason reasons as in proof in Appendix A and $\Pi^{*}(x,\eta)$ is concave (see  Lemma 1). Note that $W(x_{s},\eta)=\Pi^{*}(x_{s},\eta)$.  We define this function in order to be able to compute the $\frac{\partial \Pi^{*}(x_{s},\eta)}{\partial x}$ in terms of the derivative of $W(x,\eta)$ using the result in \cite{Benveniste} \cite{Rockafellar}.

 $\lambda_{s}=-h(x_{s})$ is in the interior so there is a small neighborhood around $-h(x_{s})$, which is also in the interior of  $[0,\lambda^{sup}]$. Let us consider another policy in which the online firm has $-h(x)$ sponsored arrivals if $ 0 \leq -h(x)\leq \lambda^{sup} $. If $-h(x) \geq \lambda^{sup}$ then the online firm gives $\lambda^{sup}$ sponsored arrivals and when $-h(x)\leq 0$ then the online firm gives $0$ sponsored arrivals. For this policy we can compute the value function in the neighborhood of $-h(x_{s})$, which is in the interior of $[0,\lambda^{sup}]$ to get $\frac{b(x)-c(-h(x))}{\rho}$ and this is equal to $W(x,\eta)$. We also know that this policy may or may not be optimal, which means $\Pi^{*}(x,\eta)\geq W(x,\eta)$ in the neighborhood of $-h(x_{s})$. 
 
 Hence, from the theorem in \cite{Benveniste} \cite{Rockafellar}  we know $\frac{\partial \Pi^{*}(x_{s},\eta)}{\partial x} = \frac{\partial W(x_{s},\eta)}{\partial x} $. We also know that $c^{'}(\lambda_{s})= \frac{\partial \Pi^{*}(x_{s},\eta)}{\partial x}$. Therefore, substituting $\lambda_{s}=-h(x_{s})$ we have
 \begin{eqnarray*}
  c^{'}(-h(x_{s}))= \frac{\partial W(x_{s},\eta)}{\partial x} \\
    c^{'}(-h(x_{s}))= \frac{b^{'}(x_{s})+c^{'}(-h(x_{s}))\frac{\partial h_{\eta}(x_{s})}{\partial x}}{\rho} \\
    \frac{b^{'}(x_{s})}{\rho - \frac{\partial h_{\eta}(x_{s})}{\partial x}} = c^{'}(-h_{\eta}(x_{s}))
  \end{eqnarray*}

The solution of the above equation gives the steady state user base level $x_{s}$. . We are interested in analyzing the  steady state value when the network effects decrease. Basically as the network effects decrease, we know that the steady state still exists (because the result in Theorem 1 is true for all values of $\eta \leq \eta^{sup}$). Next, we will derive an expression to compute the rate of change of $x_{s}$. In this proof we will make the dependence of $x_{s}$ and $\lambda_{s}$ on $\eta$ explicit for clarity and denote them as $x_{s} (\eta)$  and $\lambda_{s}(\eta)$ respectively.  Note that $x_{s}(\eta)$ is a well defined function because we know that the steady state exists and is bounded for all the values $\eta \leq \eta^{sup}$. In order to analyze the $x_{s} (\eta)$ as a function $\eta$, we first need to show that $x_{s} (\eta)$ is continuously differentiable in $\eta$.
  Consider the steady state equation 
  \begin{equation}
  \frac{b^{'}(x_{s}(\eta))}{\rho - \frac{dh_{\eta}(x_{s}(\eta))}{dx}} = c^{'}(-h(x_{s}(\eta)))
  \end{equation}
   The smoothness properties of the functions in the LHS and RHS of the steady state equation (3) can be used to show that $x_{s}(\eta)$ is continuously differentiable.
  Let us obtain the expression for $\frac{\partial x_{s}}{\partial \eta}$ by differentiating both the LHS and the RHS in the above equation (3), where $\frac{\partial }{\partial x}$ is the partial derivative w.r.t $x$.
  \begin{eqnarray*}
  && \frac{b^{''}(x_{s}(\eta))}{\rho - h_{\eta}^{'}(x_{s}(\eta))} \frac{\partial x_{s}(\eta)}{\partial \eta} + \frac{b^{'}(x_{s}(\eta))}{(\rho - h_{\eta}^{'}(x_{s}(\eta)))^2}(\frac{\partial h_{\eta}^{'}(x_{s}(\eta))}{\partial x_{s}(\eta)}\frac{\partial x_{s}(\eta)}{\partial \eta} + \frac{\partial h_{\eta}^{'}(x_{s}(\eta))}{\partial \eta})= \\
  && -c^{'}(-h_{\eta}(x_{s}(\eta)) (\frac{\partial h_{\eta}(x_{s}(\eta))}{\partial x_{s}(\eta)}\frac{dx_{s}}{d\eta} + \frac{dh_{\eta}(x_{s})}{d\eta}) 
  \end{eqnarray*}
  In the above expression we can rearrange the terms to obtain $\frac{\partial x_{s}(\eta)}{\partial \eta}$. 
  \begin{eqnarray*}
  && \frac{\partial x_{s}(\eta)}{\partial \eta} (\frac{b^{''}(x_{s}(\eta))}{\rho  - h_{\eta}^{'}(x_{s}(\eta))} +\frac{b^{'}(x_{s}(\eta))}{(\rho - h_{\eta}^{'}(x_{s}(\eta)))^2}(\frac{\partial h_{\eta}^{'}(x_{s}(\eta))}{\partial x_{s}(\eta)}) + c^{'}(-h_{\eta}(x_{s}(\eta)))\frac{\partial h_{\eta}(x_{s}(\eta))}{\partial x_{s}(\eta)}) = \\
  &&  -\frac{\partial h_{\eta}^{'}(x_{s}(\eta))}{\partial \eta}\frac{b^{'}(x_{s}(\eta))}{(\rho-h_{\eta}^{'}(x_{s}(\eta))^2)}-c^{'}(-h_{\eta}(x_{s}(\eta)))\frac{\partial h_{\eta}(x_{s}(\eta))}{\partial \eta}
  \end{eqnarray*}
  
  We know that $\lambda_{s}(\eta) =  -h_{\eta}(x_{s}(\eta))$. Hence, $\frac{\partial \lambda_{s}(\eta)}{\partial \eta} = - (\frac{\partial h_{\eta}(x_{s}(\eta))}{\partial x_{s}(\eta)} \frac{\partial x_{s}(\eta)}{\partial \eta} + \frac{\partial h_{\eta}(x_{s}(\eta))}{\partial \eta})$. Substituting the expression for $\frac{\partial x_{s}(\eta)}{\partial \eta}$ in the expression for $\frac{\partial \lambda_{s}(\eta)}{\partial \eta}$ we can obtain an expression for $\frac{\partial \lambda_{s}(\eta)}{\partial \eta}$ that depends on $x_{s}(\eta)$ and other known functions. Hence, the only unknown term is $x_{s}(\eta)$. Our objective is to find the sign of $\frac{\partial \lambda_{s}(\eta)}{\partial \eta}$. Simplifications of the above expressions yield that the sign of  $\frac{\partial \lambda_{s}(\eta)}{\partial \eta}$ is the same as 
  $b^{''}(x_{s}(\eta))\frac{\partial h_{\eta}(x)}{\partial \eta} + -\frac{b^{'}(x_{s}(\eta))}{\rho - h_{\eta}^{'}(x_{s}(\eta))}(\frac{\partial h_{\eta}^{'}(x_{s}(\eta))}{\partial \eta}\frac{\partial h_{\eta}(x_{s}(\eta))}{\partial x_{s}(\eta)} - \frac{\partial h_{\eta}^{'}(p_{s}(\eta))}{\partial x_{s}(\eta)}\frac{\partial h_{\eta}(x_{s}(\eta))}{\partial \eta} )$. If $\frac{\partial \lambda_{s}(\eta)}{\partial \eta}<0$ then $b^{''}(x_{s}(\eta))\frac{\partial h_{\eta}(x_{s}(\eta))}{\partial \eta} - \frac{b^{'}(x_{s}(\eta))}{\rho - h_{\eta}^{'}(x_{s}(\eta))}(\frac{\partial h_{\eta}^{'}(x_{s}(\eta))}{\partial \eta}\frac{\partial h_{\eta}(x_{s}(\eta))}{\partial x_{s}(\eta)} - \frac{\partial h_{\eta}^{'}(x_{s}(\eta))}{\partial x_{s}(\eta)}\frac{\partial h_{\eta}(x_{s}(\eta))}{\partial \eta} )< 0$, which implies $\frac{b^{''}(x_{s}(\eta))}{b^{'}(x_{s}(\eta))}x_{s}(\eta) <(\frac{\partial h_{\eta}^{'}(x_{s}(\eta))}{\partial \eta}\frac{\partial h_{\eta}(x_{s}(\eta))}{\partial x_{s}(\eta)} - \frac{\partial h_{\eta}^{'}(x_{s}(\eta))}{\partial x_{s}(\eta)}\frac{\partial h_{\eta}(x_{s}(\eta))}{\partial \eta} )\frac{x_{s}(\eta)}{(\rho -h_{\eta}^{'}(x_{s}(\eta)))\frac{\partial h_{\eta}(x_{s}(\eta))}{\partial \eta}} $. We can further rearrange the equation as follows
  $\rho < h_{\eta}^{'}(x_{s}(\eta)) + (\frac{\partial h_{\eta}^{'}(x_{s}(\eta))}{\partial \eta}\frac{\partial h_{\eta}(x_{s}(\eta))}{\partial x_{s}(\eta)} - \frac{\partial h_{\eta}^{'}(x_{s}(\eta))}{\partial x_{s}(\eta)}\frac{\partial h_{\eta}(x_{s}(\eta))}{\partial \eta} )\frac{1}{\frac{\partial h_{\eta}(x_{s}(\eta))}{\partial \eta}\frac{b^{'}(x_{s}(\eta))}{b^{''}(x_{s}(\eta))}} $. This condition in the previous equation is sufficient to deduce the sign of $\frac{\partial \lambda_{s}(\eta)}{\partial \eta}$ but the condition depends on $x_{s}(\eta)$, which we do not know. Define the infimum of the term on the RHS of the equation, $\inf_{x\in [0,x_{u}]} m(x)$, where $m(x)=\big(h_{\eta}^{'}(x) + (\frac{\partial h_{\eta}^{'}(x)}{\partial \eta}\frac{\partial h_{\eta}(x)}{\partial x} - \frac{\partial h_{\eta}^{'}(x)}{\partial x}\frac{\partial h_{\eta}(x)}{\partial \eta} )\frac{1}{\frac{\partial h_{\eta}(x)}{\partial \eta}\frac{b^{'}(x)}{b^{''}(x)}}\big)$. Define a threshold for $\rho_{1}^{l}= \inf_{x\in [0,x_{u}]} m(x)$ Hence, if $\rho \leq \rho_{1}^{l}$ then  $\frac{\partial \lambda_{s}(\eta)}{\partial \eta} <0$. The above result holds for the uniform threshold $\rho \leq \rho^{l}$, see Appendix E.
  
  Similarly, we can define an upper threshold for $\rho$ as well $\sup_{x\in [0,x_{u}]} m(x)$ and  $\rho_{1}^{u}= \sup_{x\in [0,x_{u}]} m(x)$. On the same lines as above it can be shown that if $\rho \geq \rho^{u}$ then $\frac{\partial \lambda_{s}(\eta)}{\partial \eta} >0$.

  Next, we prove Theorem 6. 
  We know that
  \begin{equation} 
  \frac{b^{''}(x_{s}(\eta))}{b^{'}(x_{s}(\eta))}x_{s}(\eta) <(\frac{\partial h_{\eta}^{'}(x_{s}(\eta))}{\partial \eta}\frac{\partial h_{\eta}(x_{s}(\eta))}{\partial x_{s}(\eta)} - \frac{\partial h_{\eta}^{'}(x_{s}(\eta))}{\partial x_{s}(\eta)}\frac{\partial h_{\eta}(x_{s}(\eta))}{\partial \eta} )\frac{x_{s}(\eta)}{(\rho -h_{\eta}^{'}(x_{s}(\eta)))\frac{\partial h_{\eta}(x_{s}(\eta))}{\partial \eta}} 
  \end{equation} 
  is equivalent to $\frac{\partial \lambda_{s}(\eta)}{\partial \eta} <0$. 
  
  If we take the supremum of the term in the LHS of the inequality (4) and the infimum of the term in the RHS in (4), then it will be sufficient condition to ensure that $\frac{\partial \lambda_{s}(\eta)}{\partial \eta} <0$. Define the supremum of the term in LHS as follows $\sup_{x\in [0,x_{u}]}\frac{b^{''}(x) x}{b^{'}(x)}$ and the infimum of the term in the RHS as follows $\Theta = \inf_{x \in [0,x_{u}], \eta \in [0,\eta^{sup}]}(\frac{\partial h_{\eta}^{'}(x)}{\partial \eta}\frac{\partial h_{\eta}(x)}{\partial x} - \frac{\partial h_{\eta}^{'}(x)}{\partial x}\frac{\partial h_{\eta}(x)}{\partial \eta} )\frac{x}{(\rho -h_{\eta}^{'}(x)\frac{\partial h_{\eta}(x)}{\partial \eta}}$. Therefore,  if $\sup_{x\in [0,x_{u}]}\frac{b^{''}(x) x}{b^{'}(x)} < \Theta$  then  $\frac{\partial \lambda_{s}(\eta)}{\partial \eta} <0$. 
  
  Similarly, we can define $\Delta = \sup_{x \in [0,x_{u}], \eta \in [0,\eta^{sup}]}(\frac{\partial h_{\eta}^{'}(x)}{\partial \eta}\frac{\partial h_{\eta}(x)}{\partial x} - \frac{\partial h_{\eta}^{'}(x)}{\partial x}\frac{\partial h_{\eta}(x)}{\partial \eta} )\frac{x}{(\rho -h_{\eta}^{'}(x)\frac{\partial h_{\eta}(x)}{\partial \eta}}$. On the same lines it can be shown that $\inf_{x\in [0,x_{u}]}\frac{b^{''}(x) x}{b^{'}(x)} > \Delta $ then $\frac{\partial \lambda_{s}(\eta)}{\partial \eta} >0$. This proves Theorem 6.

\section{More general benefit functions}

\subsection{Model and results overview}

In the Sections 3-6 we assumed that the benefit function  is concave. We now discuss the extension of our model to more general class of benefit functions. The benefit function $b(x)$ considered is a continuously differentiable increasing function that is not necessarily concave upto a threshold $t_{threshold}$ and beyond the threshold the benefit function is concave. Formally stated $b:[0,t_{threshold}] \rightarrow \mathbb{R}_{+}$ is convex and $b:[t_{threshold}, \infty) \rightarrow \mathbb{R}_{+}$ is concave. The benefit function is continuously differentiable and increasing on the entire domain $[0,\infty)$. 
An example of such a function is sigmoid function which is S-shaped. We will assume that the threshold $t_{threshold}$ is not very high, i.e. $t_{threshold}<f^{inf}\tilde{\eta}$. This assumption implies that the threshold is less than the steady state user base level. This is reasonable because closer to the steady state the market for the online firm starts saturating and that is why the benefits will also saturate. 

For the above class of benefit functions we can extend the following results discussed in Sections 3-6. Considering more general benefit functions do not affect the proofs of Proposition 1, Proposition 2, Theorem 1 and, Theorem 2 do not change much and these results continue to hold.
In Theorem 3 we had shown that the optimal policy of advertisements and referrals was monotonically decreasing. This does not completely extend to this scenario. The reason for this is that the marginal benefit of adding another user when the optimal user base trajectory is in the regime $[0,t_{threshold}]$ can increase due to the convexity of benefit function in the regime. However, it can be shown that once the optimal user base trajectory enters the regime $[t_{threshold}, \infty)$ the advertisements and referrals decrease with time, which is due to the concavity of the benefit function in this regime.
The comparisons in Theorem 4, which show the impact of network effects on the optimal policy as a function of user base continue to hold. Also, the comparisons of the optimal policy in the steady state under a change in network effects shown in Theorem 5 and 6 continue to hold. Next we give a detailed discussion of results.

\subsection{Detailed discussion of results}

In the model presented in Sections 3-6 we had assumed that the benefit function is concave. Consider a more general benefit function defined as $b:[0,t_{threshold}]\rightarrow \mathbb{R}_{+}$, where $b(.)$ is a continuously differentiable and increasing function (not necessarily concave) on the domain $[0,t_{threshold}]$. In the domain $[t_{threshold},\infty)$ the function $b:[t_{threshold}, \infty)\rightarrow \mathbb{R}_{+}$ is concave and increasing. On the entire domain $[0,\infty)$ the benefit function is continuously differentiable and increasing. 
Next, we discuss how do the results presented in Section 3-6 change for more general benefit functions.
In Proposition 1 we discussed the existence and uniqueness of optimal policy as a function of user base $\zeta^{*}$. It can be shown that the optimal policy as a function of time exists using the same steps as in Appendix A. The fact that benefit function is not concave does not impact the conditions required for existence in \cite{existence}. If this is true then the optimal value function $\Pi^{*}(x,\eta)$ is well defined. If the value function $\Pi^{*}(x,\eta)$ is well-defined then we can follow exactly the same steps as given in Appendix A to show that the optimal policy is unique. 
Since the optimal policy as a function of user base $\zeta^{*}$ is unique we will use this to solve for the optimal policy as a function of time. The optimal policy as a function of time will be equal to $\zeta^{*}(x_{x_{initial}}^{*}(t))$, where $x_{x_{initial}}^{*}(t)$ is the optimal user base trajectory. Next, we would need to show that the optimal user base trajectory as a function of time exists and is unique. Due to Assumption 2, we can show that the solution to the differential equation of the user base exists and is unique. This will establish the uniqueness of the optimal policy as a function of time. 

We next discuss how do the theorems change in this setting.

In Theorem 1 we showed the existence and uniqueness of the steady state. In order to prove the existence we  had exploited the property that the user base trajectory is bounded and increasing. In this case as well the user base trajectory is bounded above and increasing for exactly the same reasons as mentioned in Appendix A, therefore, steady state will exist. In order to show uniqueness we will exploit the concavity of the $\Pi^{*}(x,\eta)$ in the neighborhood of the steady state (justified in next paragraph) and use the same steps as before.
The same idea can also be used to show the convergence as in Theorem 2 (increasing and bounded function will converge). 

In Theorem 3 we had shown that the optimal policy was a decreasing function of time. For this we relied on showing that $\Pi^{*}(x,\eta)$ is concave in $x$. In this case however, it is not necessarily true that $\Pi^{*}(x,\eta)$ is concave on the entire domain $[0,\infty)$. This is because the benefit function is not concave in the domain $[0,t_{threshold}]$. However, since  the benefit function is concave in the regime $[t_{threshold},\infty)$ we can show that $\Pi^{*}(x,\eta)$ is concave in the regime $[t_{threshold},\infty)$. In order to show this we can write the net benefit function $b(x)-c(\frac{dx}{dt}-h_{\eta}(x))$, which is concave in $[t_{threshold},\infty)$. Note that once the user base trajectory crosses the $t_{threshold}$ it will stay above $t_{threshold}$, this is because $h_{\eta}(t_{threshold})>0$ (follows from the assumption that $t_{threshold} < f^{inf} \tilde{\eta} $). Since the user base trajectory always stays in the regime in which the net benefit function is concave and since  the regime itself is a convex set, we can use the same result as in \cite{Benveniste} to show that the $\Pi^{*}(x,\eta)$ is concave in $[t_{threshold},\infty)$. Also, note that we can show that the user base trajectory is increasing, this is because the rate of change of user base $h_{\eta}(x) + \zeta^{*}(x)>0$ at all points before the steady state. The steady state value has to be greater than $f^{inf}\tilde{\eta} > t_{threshold}$, which means that the user base trajectory will at some point exceed $t_{threshold}$.
 Since $\Pi^{*}(x,\eta)$ is concave in 
$[t_{threshold},\infty)$ and this combined with the fact that the user base trajectory is increasing and will cross the threshold $t_{threshold}$ we can conclude that the policy will decrease with time after some time (when user base exceeds $t_{threshold}$).

In Theorem  4  we show the impact of the level of foresight and the size of the user base on the online firm's decision when it sees a decline in the network effects. In Appendix E, where we show the proof for Theorem 4  we only rely on the fact that the benefit function is increasing. So Theorem 4 continue to hold in this setting.  

In Theorems 5  and 6 we analyze the impact of change in network effects on the optimal policy  as a function of time. In the proof shown in Appendix F we developed the steady state equation. In developing the steady state equation we used the fact the value function $\Pi^{*}(x,\eta)$ is concave in the neighborhood around the steady state value. Note that the steady state value is more than $t_{threshold}$. Since $\Pi^{*}(x,\eta)$ is concave in the region $[t_{threshold}, \infty )$, this means that $\Pi^{*}(x,\eta)$ is concave in a neighborhood around steady state. Therefore, we can develop the same steady state equation and the results analogous to Theorem 6 and 7 can be arrived at.  

\section{Users take decisions to join the online firm}
\subsection{Model and results overview}

In this section we extend our model to allow for users decisions.  We first describe the additions in the model discussed in Section 3. 
In Section 3 we had assumed that the direct arrival rate was $f(x)$. We need to modify the direct arrivals to take users' decisions into account. Suppose $f(x)$ users hear about the online firm (for e.g. word of mouth) per unit time and have to decide whether to enter the online firm or not \footnote{$f(x)$ users can also comprise of the users who already know about the online firm and need to decide whether to visit the online firm or not}. Each user is rational and has a preference $\alpha$ for the online firm drawn from a uniform distribution $[\alpha_{min},\alpha_{max}]$ \footnote{Uniform distribution is not a restrictive assumption and the results discussed in this subsection extend to more general distributions as well}. The user's benefit from joining the online firm increases with the size of the user base. The user's benefit also increases with the preference of the user for the online firm and it is given as $\alpha x$. Each user bears a fixed cost for joining the online firm given as $c_{enter}$. The utility of the user from entering the online firm is given by $  \alpha x -c_{enter} $ and is zero otherwise. At every time instance  $f(x)$ users need to decide whether to visit the online firm and suppose a fraction $w(x)$ of these users decide to join the online firm per unit time, then the direct arrival rate is given by $f(x)w(x)$. 

In Section 3 we had the sponsored arrival rate at time $t$ given as $\lambda(t)$, which was a result of the advertisements and referrals given by the online firm.  Suppose the online firm can control the number of users who come across the advertisements $\lambda(t)$ users. These users need to decide whether to enter the online firm or not. Hence, if a fraction $w(x)$  of users decide to enter the online firm the rate of the sponsored arrivals is $\lambda(t).w(x)$.  Therefore, the rate at which the user base of the online firm changes is given as 
$\frac{dx}{dt} = f(x)w(x) + \lambda(t)w(x) - \frac{1}{\tilde{\eta} +\eta g(x)}x $. 
The benefits and costs for the online firm remain the same as in the Section 3. Hence, the long-term discounted average profit of the online firm is given as $B(\bar{x})-C(\bar{\lambda})$ (as defined in Section 3). 

Before we discuss the equilibrium of the above model we make a remark about the model discussed above. In the model discussed above we can also include the payments made to the online firm directly in the form of subscription fees. We can do so by adding a constant fee to the cost of entry for the users. We also need to modify the benefit function of the online firm as follows. The benefit function will be sum of $b(x)$, which represents the benefits generated by users clicking on advertisements, and a term that is equal to the subscription fees times the rate of arrivals (sum of the direct and the sponsored). We do not add the subscription costs to avoid a complicated description. 

\subsubsection{Equilibrium analysis:}

In the model discussed above both the users and the online firm make decisions and want to maximize their own utilities. We want to compute a policy for the online firm and the users in which both  the online firm and the users cannot deviate unilaterally to increase their utilities. 

Let us consider a user with a preference level $\alpha$ for the online firm. Suppose the user comes across the online firm directly or through an advertisement then he needs to decide whether to enter the online firm or not. The user has two actions $\{\text{enter}, \text{not enter}\}$. We assume that the user knows the size of the user base $p$.   If  the utility from entering is greater than zero $\alpha.x -c_{enter} > 0$  then the unique best response of the user is to enter the online firm and not enter other wise. Note that the user is assumed to be myopic and hence, the best response only depends on the size of the user base $p$ and not on the future decisions to be made by the online firm. Hence, the best response of the user as discussed above gives the strategy of the user in equilibrium. 
Let the set of $\alpha$'s which satisfy $\alpha.x -c_{enter} > 0$ be given as $[\alpha(x),\alpha_{max}]$. Hence, the fraction of users who enter the online firm is given as $w^{eq}(x)=\frac{(\alpha_{max}-\alpha(x))}{\alpha_{max}-\alpha_{min}}$. 

The online firm has a belief about each user's actions.  The online firm also knows the preference distribution of the users.  We assume that the online firm holds the same belief for users of the same type $\alpha$. Hence, the online firm will have a belief about the fraction of the users that will enter, denoted as $\hat{w}(x)$. The online firm's profit maximization problem for the belief function $\hat{w}(x)$ is stated as follows. 

\begin{eqnarray*}
	&&  \hat{\Pi}(x_{initial}, \eta) = \max_{0 \leq \lambda(t) \leq \lambda^{sup}, \forall t \geq 0,\;} B(\bar{x})-C(\bar{\lambda}) \\
	&& \text{subject to}\; \frac{dx}{dt} = f(x)\hat{w}(x) + \lambda(t)\hat{w}(x) - \frac{1}{(\tilde{\eta} +\eta g(x))}x,\\
\end{eqnarray*}

The optimal solution of the above problem is denoted as $\hat{\lambda}_{x_{initial}}(t ,\hat{w}(.))$ and the corresponding user base trajectory is given as $\hat{x}_{x_{initial}}(t ,\hat{w}(.))$. Note that the optimal policy is a function of the belief $\hat{w}(.)$. In equilibrium the fraction of users joining the online firm is given as $w^{eq}(x)$. The online firm computes the optimal policy in equilibrium with the belief $\hat{w}(x)=w^{eq}(x)$. We already know that the users will not want to deviate from their decisions that result in $w^{eq}(x)$ fraction of users entering. The online firm will not want to deviate either as it maximizes its profit. The equilibrium policy as a function of time is denoted as $\lambda^{eq}_{x_{initial}}$ and the corresponding user base dynamic is denoted as $x^{eq}_{x_{initial}}$ 

Next, we can also formulate the HJB equation for the above problem in the same manner as in Section 3. Let us denote the optimal policy as a function of  user base as $\zeta^{eq}_{x_{initial}}(.)$.

We can show that the results presented in Section 4-6 continue to hold under certain assumptions. We now give the intuition of how to reduce the above analysis for the case above to the analysis in Section 4-6. 
Let us consider the differential equation in Section 3 and the one developed above. The first thing to note is that the function $f(x)$ represented the direct arrivals  in Section 3. Here, we have a new function which is $f(x)w(x)$.  Earlier we did not have a term multiplied with $\lambda(t)$, while here we have a term of $w(x)$ multiplied with $\lambda(t)$. This is an important difference in comparison to the earlier sections. We can define a function $\hat{h}_{\eta}(x) = f(x) - \frac{1}{w(x)(\bar{\eta}+\eta g(x))}x$. If we replace $h_{\eta}(x)$ in Sections 4-6 with  $\hat{h}_{\eta}(x)$ and also add some more regularity conditions then the main results presented continue to hold. Next, we give a detailed discussion about the results.

\subsection{Detailed discussion of results}

As pointed out earlier in this case there are two differences in this extension. First, the rate of direct arrivals changes from $f(x)$ to $f(x)w(x)$. $w(x)$ is increasing and bounded above, which is implied from its definition. We can assume in addition that $w(x)$ is continuously differentiable and also assume that it is bounded below by a positive constant. We assume that $\alpha_{min}=0$ and $\alpha_{max}=1$ without loss of generality. If the distribution of the users is uniform then a user enters if $\alpha \geq \frac{c_{enter}}{x}$.  Hence, if $x\geq c_{enter}$ then only will a user enter. The fraction of users entering if $x \geq c_{enter}$ is $w(x)= (1-\frac{c_{enter}}{x})$. Now if $x_{initial}> c_{enter}(1+\kappa)$ then the function $w(x)$ is continuously differentiable and is lower bounded by $ \frac{\kappa}{1+\kappa}$. 
Second, the rate of sponsored arrivals changes from $\lambda(t)$ to $\lambda(t)f(x)$

In Proposition 1 and 2 we showed the existence and uniqueness of the optimal policy. First we showed that the optimal policy as a function of time exists. For that we used the sufficient conditions shown in \cite{existence}.  Even under the two key changes described above the sufficient conditions in  \cite{existence} continue to hold.
  Therefore, the optimal policy as a function of time exists and  $\Pi^{*}(x,\eta)$ is well defined. We use this to show the existence and uniqueness of the optimal policy as a function of user base $\zeta^{*}(.)$. We need to modify the part 1 of Assumption 1 as follows $c^{'}(0)< \frac{d\Pi(x,\eta)}{dx} w(x) < c^{'}(\lambda^{sup})$. This modification ensures that the optimal policy is in the interior of $[0,\lambda^{sup}]$. $\Pi^{*}(x,\eta)$ has to satisfy the HJB equation and the optimal policy as a function of user base will be the maximizer of the RHS of the HJB equation given as $c'^{-1}(\frac{d\Pi(x,\eta)}{dx} w(x))$. This proves Proposition 1.

Since we already showed the existence of the optimal policy as a function of time, we need to show uniqueness next. Basically we need to show that $\frac{dx}{dt} = f(x)w(x) +  c'^{-1}(\frac{d\Pi(x,\eta)}{dx} w(x)) w(x) - \frac{1}{\tilde{\eta}+\eta g(x)}x$ for initial value $x_{initial}$ has a unique solution. We need to modify assumption 2 as follows. Define $\hat{h}_{\eta}(x) = f(x)w(x)  -\frac{1}{\tilde{\eta}+\eta g(x)} $ and the assume that $\hat{h}_{\eta}(x)$ is continuously differentiable. $c'^{-1}(\frac{d\Pi(x,\eta)}{dx} w(x)) w(x)$ can be argued to be continuously differentiable based on $w(x)$ and $\frac{d\Pi(x,\eta)}{dx}$ continuous differentiability. Thus we can use the results from theory of ODEs (in the same way as we did in Appendix A) \cite{Levinson} to show the uniqueness. 

In Theorem 1 we showed the existence of the steady state. Note that in this case as well we can show that the user base trajectory is bounded above. This is because the $f(x)w(x)$ is bounded above by $f^{sup}$ and $c'^{-1}(\frac{d\Pi(x,\eta)}{dx} w(x)) w(x)$ is bounded above by $\lambda^{sup}$. Therefore, the rate of change of user base is bounded by $f^{sup}+\lambda^{sup} - \frac{1}{\tilde{\eta}+\eta g(x)}x$. Every trajectory that has a rate of change given by $f^{sup}+\lambda^{sup} - \frac{1}{\tilde{\eta}+\eta g(x)}x$ is bounded. Recall that we had the assumption that  $x_{initial}<f^{inf}\tilde{\eta}$, which results in $h_{\eta}(x_{initial})>0$. This assumption was used in proving that the steady state exists.
We can modify this assumption to $x_{initial}<f^{inf}\frac{\kappa}{1+\kappa}$. Under this modified assumption the rate of change of user base $\hat{h}_{\eta}(x_{initial}) >0$. In order to show that the steady state exists we need to show that  $f(x)w(x) +  c'^{-1}(\frac{d\Pi(p,\eta)}{dx} w(x)) w(x) - \frac{1}{\tilde{\eta}+\eta g(x)}x=0$ has a solution. We already know that $f(x)w(x) +  c'^{-1}(\frac{d\Pi(x,\eta)}{dx} w(x)) w(x) - \frac{1}{\tilde{\eta}+\eta g(x)}x>0$ when $x=x_{initial}$. $f(x)w(x) +  c'^{-1}(\frac{d\Pi(x,\eta)}{dx} w(x)) w(x) - \frac{1}{\tilde{\eta}+\eta g(x)}x$ at $x=x_{u}$ has to be negative. Because $f(x)w(x) +  c'^{-1}(\frac{d\Pi(p,\eta)}{dx} w(x)) w(x) - \frac{1}{\tilde{\eta}+\eta g(x)}x$ is continuous it will have a root in the interval $[x_{initial}, x_{u}]$. This proves the existence of the steady state. 

In Theorem 2 we showed the convergence of the policy to the steady state. We can replicate the same steps as given in Appendix C and prove convergence

In Theorem 3 we showed that the optimal policy will be decreasing as a function of time. Here, it is not necessary that the optimal policy will decrease with time. This is because the optimal policy as a function of user base is proportional to $\frac{\partial \Pi^{*}(x,\eta)}{\partial x}w(x)$ and $\frac{\partial \Pi^{*}(x,\eta)}{\partial x}w(x)$ does not necessarily decrease with user base. We will show that $\frac{\partial \Pi^{*}(x,\eta)}{\partial x}$ is concave, which will imply that $\frac{\partial \Pi^{*}(x,\eta)}{\partial x}$ is decreasing. But $w(x)$ is increasing. Therefore, if the percentage decrease in $\frac{\partial \Pi^{*}(x,\eta)}{\partial x}$ is greater than the percentage increase in $w(x)$ then the optimal policy decreases otherwise it increases.

Next, we discuss the impact on Theorem 4. We can write the HJB equation for this case as follows. 
\begin{equation*}
\rho \Pi^{*}(x,\eta) = b(x) - c(q(\frac{\partial \Pi^{*}(x,\eta)}{\partial x}w(x))) + \frac{\partial \Pi^{*}(x,\eta)}{\partial x} w(x )(f(x)+ q(\frac{\partial \Pi^{*}(x,\eta)}{\partial x}w(x)) - \frac{1}{w(x)(\tilde{\eta}+\eta g(x))}x)   
\end{equation*}

As pointed out above that the optimal policy is proportional to  $y(x,\eta) =\frac{\partial \Pi^{*}(x,\eta)}{\partial x}w(x)$. We take the derivative of the above HJB equation w.r.t $\eta$ to obtain the following expression. 

\begin{eqnarray*}
&& \rho \frac{\partial \Pi^{*}(x,\eta)}{\partial \eta} = -c'(q'(y(x)))\frac{\partial y(x)}{\partial \eta} + \frac{\partial y(x)}{\partial \eta}(f(x)+ q(y(x)) - \frac{1}{w(x)(\tilde{\eta}+\eta g(x))}x)   + y(x)q^{'}(y(x))\frac{\partial y(x)}{\partial \eta}  \\
&& + y(x)\frac{g(x)x}{w(x)(\tilde{\eta}+ \eta g(x))^2}
\end{eqnarray*}

Rearranging the above equation we get the following
\begin{eqnarray*}
&& \rho \frac{\partial \Pi^{*}(x,\eta)}{\partial \eta} -y(x)\frac{g(x)x}{w(x)(\tilde{\eta}+ \eta g(x))^2} = \frac{\partial y(x)}{\partial \eta}( -c'(q'(y(x))) + f(x)+ q(y(x)) - \frac{1}{w(x)(\tilde{\eta}+\eta g(x))}x  \\
&&  + y(x)q^{'}(y(x)))
\end{eqnarray*}
Now we can carry out an analysis on the same lines as shown in Appendix E. Basically,  a sufficiently low and high discount rate $\rho$ will lead to a negative or a positive LHS. 
The sign of the term $( -c'(q'(y(x))) + f(x)+ q(y(x)) - \frac{1}{w(x)(\tilde{\eta}+\eta g(x))}x  + y(x)q^{'}(y(x)))$ will depend on the size of the user base. 
This way we can prove Theorem 4. 

Next, we discuss how Theorems 5 and 6 change. For Theorems 5 and 6 we had arrived at the steady state equations using the concavity of the value function. In this case as well we can show that the value function is concave in a manner similar to Lemma 1. We need to assume that $b(x)-c(\frac{dx}{dt}\frac{1}{w(x)}-f(x)+ \frac{1}{(\tilde{\eta}+\eta g(x))w(x)}x)$ is jointly concave in $(x,\frac{dx}{dt})$ where $x$ is restricted to be in $[x_{initial}, x_{u}]$. Using the sufficient conditions for concavity in \cite{Benveniste} we thus have $\Pi^{*}(x,\eta)$ to be concave. 
We can therefore derive the steady state equation for this case as follows. 
\begin{eqnarray*}
	c^{'}(-\hat{h}_{\eta}(x_{s}))\frac{1}{w(x_{s})} =  \frac{b^{'}(x_{s})+c^{'}(\hat{h}_{\eta}(x_{s}))\frac{\partial \hat{h}_{\eta}(x_{s})}{\partial x}}{\rho}
\end{eqnarray*}

In the above $\hat{h}_{\eta}(x_{s}) = f(x) - \frac{1}{(\tilde{\eta} + \eta g(x))w(x)} x$. The above equation can be analyzed in the similar manner as in Appendix F. The analysis will show that the heterogeneity in the benefit function determines the impact of network effects. 

\section{Stochastic arrival/exits}
In this section we allow for uncertainty in the arrivals and exits of users. This modifies the user user base differential equation to the following stochastic differential equation, $dX_{t} = (\theta +\lambda(t)-\lambda_{d}X)dt+\sigma X dW_{t}$. ($\lambda_{d}=\frac{1}\eta$). The HJB equation corresponding to this case is
\begin{eqnarray*}
	\rho \Pi(x) = \max_{\lambda}(\Gamma^2-(\Gamma-x)^2-c\lambda^2+\Pi^{'}(x)(\theta+\lambda-\lambda_{d}x) + \frac{1}{2}\Pi^{''}(x)\sigma^2 )
\end{eqnarray*}
In the above setting we allow for negative referrals for analytical tractability, but we do not allow that the user base becomes negative. Hence, the optimal $\lambda$ should ensure that the user base is always positive. Solving the optimal $\lambda$ we get $\lambda=\frac{\Pi^{'}(x)}{2c}$.
Let's assume that $\Pi(x)=A^{'}x^2+B^{'}x+C^{'}$ and substitute $\lambda=\frac{\Pi^{'}(x)}{2c}=\frac{2A^{'}x+B^{'}}{2c}$ in the above HJB equation. Equating coefficients on both the sides we get $A^{'} = \frac{2\lambda_{d}+\rho -\sigma^2-\sqrt{(2\lambda_{d}+\rho-\sigma^2)^2+4(\frac{1}{c})}}{2(\frac{1}{c})}$ and $B^{'} = 2\frac{(A^{'}\theta+\Gamma)}{(\rho+\lambda_{d}-\frac{A^{'}}{c})}$. Note that $\frac{A^{'}}{c}-\lambda_{d}<0$ and $\theta+\frac{B^{'}}{2c}>0$, these conditions together with the fact that the initial user base $x(0)\geq0$ ensure that the user base dynamic stays above zero (see \cite{zhao2009inhomogeneous}). The expected user base at a time $t$ is given as $E[x(t)] = \frac{\theta+\frac{B^{'}}{2c}}{\lambda_{d}-\frac{A^{'}}{c}}(1-e^{-(\lambda_{d}-\frac{A^{'}}{c})t}) + x(0)e^{-(\lambda_{d}-\frac{A^{'}}{c})t}$.
The expected number of referrals are given as $E[\lambda(t)] = \frac{A^{'}E[x(t)]+B^{'}}{2c}$.


If $x(0)\leq \frac{\theta+\frac{B^{'}}{2c}}{\lambda_{d}-\frac{A^{'}}{c}}$ then the expected user base $E[x(t)]$ dynamic increases with time and converges to attain $ \frac{\theta+\frac{B^{'}}{2c}}{\lambda_{d}-\frac{A^{'}}{c}}$. The corresponding trajectory of expected number of referrals will decrease with time, this is because the referrals decrease with increase in user base

We know that $\lambda = \frac{A^{'}}{c}X+\frac{B^{'}}{2c}$. We know that $B^{'}=2\frac{A^{'}\theta+\gamma}{\rho+\lambda_{d}-\frac{A^{'}}{c}}$ and from the expression observe that if $A^{'}$ increases $B^{'}$ increases as well. 

We next analyze the impact of change in network effects on the optimal policy. We can compute $\frac{\partial A^{'}}{\partial \lambda_{d}} =  c.(1- \frac{2\lambda_{d}+ \rho -\sigma^2}{\sqrt{(2\lambda_{d}+\rho-\sigma^2)^2+4\frac{1}{c}}})$. Observe that $\frac{\partial A^{'}}{\partial \lambda_{d}}\in [0,1]$. We can also obtain the expression for
 $\frac{\partial B^{'} }{\partial \lambda_{d} } = (\rho+\lambda_{d})\theta\frac{\partial A^{'}}{\partial \lambda_{d}}-A^{'}\theta + \Gamma(-1+\frac{\partial A^{'}}{\partial \lambda_{d}})$. Observe that the first term $(\rho+\lambda_{d})\theta \frac{\partial A^{'}}{\partial \lambda_{d}} >0$  and  the second term  $ \Gamma(-1+\frac{\partial A^{'}}{\partial \lambda_{d} })<0$. If $\Gamma$ is sufficiently high then $\frac{\partial B^{'}}{\partial \lambda_{d}} <0$. If $\Gamma$ is sufficiently low then $\frac{\partial B^{'}}{\partial \lambda_{d}} <0$. Hence, we can see that if $\Gamma$ is sufficiently low then both $A^{'}$ and $B^{'}$ increase with an increase in $\lambda_{d}$. Therefore, for a sufficiently low $\Gamma$ the optimal policy as a  function of the size of the user base given as $\frac{A^{'}}{c}X+\frac{B^{'}}{2c}$ is an increasing function of $\lambda_{d}$. If $\Gamma$ is sufficiently large then for sufficiently small size of the user base the optimal policy $\frac{A^{'}}{c}X+\frac{B^{'}}{2c}$ decreases, while for sufficiently large sizes of the user base the optimal policy $\frac{A^{'}}{c}X+\frac{B^{'}}{2c}$  increases.  
 
 We also want to compare the optimal policy in the steady state.  Let us denote the expected population level in the steady state as $X_{s}$. Therefore the optimal policy in the steady state is given as $\frac{A^{'}}{c}X_{s} + \frac{B^{'}}{2c}$. We now see how  $\frac{A^{'}}{c}X_{s} + \frac{B^{'}}{2c}$ changes as $\lambda_{d}$ increases. $\frac{\partial (\frac{A^{'}}{c}X_{s} + \frac{B^{'}}{2c})}{\partial \lambda_{d}} = \frac{\partial A^{'}}{\partial \lambda_{d}}  X_{s} + A^{
 	'}\frac{\partial X_{s}}{\partial \lambda_{d}}+\frac{\partial B^{'}}{\partial \lambda_{d}} $
$\frac{\partial A^{'}}{\partial \lambda_{d}}  X_{s}$. We know that $\frac{\partial A^{'}}{\partial \lambda_{d}}  X_{s}>0$, $A^{
	'}\frac{\partial X_{s}}{\partial \lambda_{d}}>0$ (because $A^{'}<0$ and $\frac{\partial X_{s}}{\partial \lambda_{d}}<0$) and $\frac{\partial B^{'}}{\partial \lambda_{d}}>0$. Therefore, increasing $\lambda_{d}$ will increase the advertisements and referrals in steady state.

We prove Theorem 7 and 8 next. 
Next, we show that as $\sigma$ increases then $A^{'}$ decreases.
\begin{eqnarray*}
	\frac{dA^{'}}{d\sigma^2} = c(-1+ \frac{(2\lambda_{d}+\rho-\sigma^2)}{\sqrt{(2\lambda_{d}+\rho-\sigma^2)^2+4(\frac{1}{c})}})
\end{eqnarray*}

Observe that $\frac{dA^{'}}{d\sigma} <0$. We know that $B^{'}$ increases with an increase in $A^{'}$.  Therefore, $\frac{dB^{'}}{d\sigma}<0$. Therefore, with increase in $\sigma$ the $\frac{2A^{'}X+B^{'}}{c}$ will be a decreasing function of $\sigma$. 


For Theorem 8 we need to compare the steady state user base level, which is $\frac{\theta+\frac{B^{'}}{2c}}{\lambda_{d}-\frac{A^{'}}{c}}$. This expression is positive and decreases with $\sigma$ as $A^{'}$ and $B^{'}$ decrease with $\sigma$. 
On the same lines we can compare the value functions for the result on profits. We obtain $C^{'}$ in the same way as $A^{'}$,$B^{'}$ and it is given as $C^{'}=B^{'}\theta +\frac{B'^{2}}{4c}$.  $C^{'}$ decreases with $\sigma$ as  $B^{'}$ decreases with $\sigma$. This combined with the definition of the value function gives the result.

\section{Concavity of the benefit function}

In this section we will show that if the users arrive in the decreasing order of their preference then the benefit generated for the online firm is a concave function that depends on the size of the user base. 

Let the potential user base be a set $S^{'}=[0,\infty)$, where each member of the set corresponds to a potential user. Also, assume that length of any interval of $S^{'}$ will denote the size of the user base that lies in that interval. Let $s$ be a potential user with a preference for the online firm $pr(s)$, where we assume without loss of generality that $pr(s)$ is  a decreasing function in $s$. We also assume that the benefit that user $s$ will generate is proportional to $pr(s)$. Let us construct another preference function that samples $pr(s)$ and looks like a staircase, which is defined as $\hat{pr}(s)$. Divide the domain into intervals of length $\Delta s$ and  $\hat{pr}(s)=pr(s_{l})$, where $s_{l}$ corresponds to the left end point of the interval to which $s$ belongs to. 
In the analysis that follows we will assume that preference of the users follows $\hat{pr}(s)$. We will show that our result will hold for every $\Delta s$. We know that $pr(s)=\lim_{\Delta s \rightarrow 0}\hat{pr}(s)$, which will extend our result to $pr(s)$. 

Consider the interval $[0,\Delta s]$. We know that the user arrivals happen in the order of the preference. All the users in the interval $[0,\Delta s]$ have the same preference $\hat{pr}(0)$ and $\hat{pr}(0)$ is the highest preference. It is clear that initially users from interval $[0,\Delta s]$ will only enter the online firm. Let the population of the users with preference $\hat{pr}(0)$ at time $t$ be given as $p_{1}(t)$, where we use subscript 1 with population to denote that these users correspond to the first interval $[0,\Delta s]$. The rate at which the population will change is given as $\frac{dx_{1}}{dt} =  f(x_{1}) + \lambda(t) - \frac{1}{\tilde{\eta}+\eta g(x_{1})}x_{1}$ as long as the population is less than $\Delta s$. 
The benefit generated by these users is given as $\hat{b}(p)= \int_{0}^{x} \hat{pr}(s) ds$, where $p$ is the size of the user base that has entered and is less than $\Delta s$. 
 If $\Delta s$ is sufficiently small, then the population trajectory of $x_{1}(t)$ will achieve $\Delta s$ level in a certain finite time. When $x_{1}(t)$ achieves $\Delta s$ then the users will still leave the user base at the rate $\frac{1}{\tilde{\eta}+\eta g(\Delta s)} \Delta s$. The total direct arrival rate is $f(\Delta s)$ and a fraction of direct arrivals will nullify the exits of the users $\frac{1}{\tilde{\eta}+\eta g(\Delta s)} \Delta s$ and the remaining direct arrivals will bring in the users from the interval $[\Delta s , 2 \Delta s]$. Note that the users from the next interval only begin entering after the users from the previous interval have already entered. Also, the size of the user base from the interval $[0,\Delta s]$ remains fixed to $\Delta s$. The users from the interval  $[\Delta s , 2 \Delta s]$ will enter at the rate $\frac{dx_{2}}{dt} =  f(x_{2}+\Delta s) + \lambda(t) - \frac{1}{\tilde{\eta}+\eta g(x_{2}+\Delta s)}(\Delta s +x_{2})$ and $\frac{dx_{1}}{dt}=0$. The total rate of change of population $\frac{dx}{dt} = f(x_{2}+\Delta s) + \lambda(t) - \frac{1}{\tilde{\eta}+\eta g(x_{2}+\Delta s)}(\Delta s +x_{2}) = f(x) + \lambda(t) - \frac{1}{\tilde{\eta}+\eta g(x)}(x) $, where $p$ is the total size of the user base. The benefit that is generated is given as $\hat{b}(x)= \int_{0}^{x} \hat{pr}(s) ds$. This argument can be generalized to all the intervals. Note that $\hat{b}(x)$ is a concave function because $\hat{pr}(s)$ is a decreasing function. We also know that the above argument holds for any $\Delta s$. Therefore, we can say that $b(x)= \lim _{\Delta s \rightarrow 0}\int_{0}^{x} \hat{pr}(s) ds = \int_{0}^{x} \lim _{\Delta s \rightarrow 0} \hat{pr}(s) ds = \int_{0}^{x} pr(s) ds$. The second equality in the above follows from monotone convergence theorem. Therefore, $b(x)$ depends only on $x$ and is concave.

 \end{APPENDICES}






%

\vspace{38 em}

\bibliographystyle{ormsv080} 
\bibliography{Management-Science-draft-Aug20}

\end{document}